\documentclass[12pt]{article}
\usepackage{theorem}
\usepackage{amssymb,amsmath,mathrsfs}
\newtheorem{prop}{}[section]

{\theorembodyfont{\upshape} }
\newcommand{\boma}[1]{{\mbox{\boldmath $#1$} }}
\begin{document}
\def\ppsi{\vartheta}
\def\PPsi{\Theta}
\def\L{M}
\def\M{N}
\def\N{V}
\def\R{R}
\def\half{{1 \over 2}}
\def\II{M}
\def\JJ{N}
\def\ga{\gamma_{\scriptscriptstyle{E}}}
\def\mmu{\nu}
\def\nnu{\mu}
\def\F{G}
\def\f{g}
\def\Fs{\mathscr{G}}
\def\Kap{\mathscr{K}}
\def\STB{\scriptsize B}
\def\STBB{\scriptsize BB}
\def\STF{\scriptsize F}
\def\STFF{\scriptsize FF}
\def\TB{\scriptsize(B)}
\def\TBB{\scriptsize(BB)}
\def\TF{\scriptsize(F)}
\def\TFF{\scriptsize(FF)}
\def\p{\alpha}
\def\pp{\beta}
\def\XXX{\mathscr X}
\def\Piu{\mathscr P}
\def\Men{\mathscr N}
\def\ffi{\varphi}
\def\ES{{\mathcal S}}
\def\J{{\mathscr J}}
\def\H{{\mathscr H}}
\def\K{{\mathscr K}}
\def\Kp{{\mathscr K}'}
\def\kp{k'}
\def\scrscr{\scriptscriptstyle}
\def\scr{\scriptstyle}
\def\dd{\displaystyle}
\def\z{w}
\def\x{w}
\def\w{\eta}
\def\B{ B_{\mbox{\scriptsize{\textbf{C}}}} }
\def\Bc{ \overline{B}_{\mbox{\scriptsize{\textbf{C}}}} }
\def\ppartial{\overline{\partial}}
\def\d{\hat{d}}
\def\TT{T}
\def\G{ {\textbf G} }
\def\Hinf{ H^{\infty}(\reali^d, \complessi) }
\def\Hn{ H^{n}(\reali^d, \complessi) }
\def\Hm{ H^{m}(\reali^d, \complessi) }
\def\Ha{ H^{\d}(\reali^d, \complessi) }
\def\Ld{L^{2}(\reali^d, \complessi)}
\def\Lpi{L^{p}(\reali^d, \complessi)}
\def\Lq{L^{q}(\reali^d, \complessi)}
\def\Lr{L^{r}(\reali^d, \complessi)}
\def\Knb{K^{best}_n}
\def\k{\mbox{{\tt k}}}
\def\D{\mbox{{\tt D}}}
\def\g{ {\textbf g} }
\def\QQQ{ {\textbf Q} }
\def\AAA{ {\textbf A} }
\def\gr{\mbox{graph}~}
\def\Q{$\mbox{Q}_a$~}
\def\PZ{$\mbox{P}^{0}_a$~}
\def\PZAL{$\mbox{P}^{0}_\alpha$~}
\def\PL{$\mbox{P}^{1/2}_a$~}
\def\PU{$\mbox{P}^{1}_a$~}
\def\PK{$\mbox{P}^{k}_a$~}
\def\PKU{$\mbox{P}^{k+1}_a$~}
\def\PI{$\mbox{P}^{i}_a$~}
\def\Pell{$\mbox{P}^{\ell}_a$~}
\def\PTM{$\mbox{P}^{3/2}_a$~}
\def\AZ{$\mbox{A}^{0}_r$~}
\def\AU{$\mbox{A}^{1}$~}
\def\epsilona{\epsilon^{\scriptscriptstyle{<}}}
\def\epsilonb{\epsilon^{\scriptscriptstyle{>}}}
\def\lgraffa{ \mbox{\Large $\{$ } \hskip -0.2cm}
\def\rgraffa{ \mbox{\Large $\}$ } }
\def\restriction{ \stackrel{\setminus}{~}\!\!\!\!|~}
\def\m{m}
\def\Fre{Fr\'echet~}
\def\I{{\mathcal N}}
\def\ap{{\scriptscriptstyle{ap}}}
\def\fiap{\varphi_{\ap}}
\def\BBB{ {\textbf B} }
\def\EEE{ {\textbf E} }
\def\FFF{ {\textbf F} }
\def\TTT{ {\textbf T} }
\def\KKK{ {\textbf K} }
\def\FFi{ {\bf \Phi} }
\def\GGam{ {\bf \Gamma} }
\def\a{a}
\def\ep{\epsilon}
\def\parn{\par\noindent}
\def\teta{M}
\def\elle{L}
\def\ro{\rho}
\def\al{\alpha}
\def\si{\sigma}
\def\be{\beta}
\def\de{\delta}
\def\la{\lambda}
\def\te{\vartheta}
\def\ch{\chi}
\def\complessi{{\textbf C}}
\def\reali{{\textbf R}}
\def\interi{{\textbf Z}}
\def\naturali{{\textbf N}}
\def\bT{{\textbf T}}
\def\T1{{\textbf T}^{1}}
\def\EE{{\mathcal E}}
\def\FF{{\mathcal F}}
\def\EFFE{{\mathscr F}}
\def\GG{{\mathcal C}}
\def\PP{{\mathcal P}}
\def\QQ{{\mathcal Q}}
\def\Np{{\hat{N}}}
\def\Lp{{\hat{L}}}
\def\Jp{{\hat{J}}}
\def\Pp{{\hat{P}}}
\def\Pip{{\hat{\Pi}}}
\def\Vp{{\hat{V}}}
\def\Ep{{\hat{E}}}
\def\Fp{{\hat{F}}}
\def\Gp{{\hat{G}}}
\def\Ip{{\hat{I}}}
\def\Tp{{\hat{T}}}
\def\Mp{{\hat{M}}}
\def\La{\Lambda}
\def\Ga{\Gamma}
\def\Si{\Sigma}
\def\Upsi{\Upsilon}
\def\Gag{{\check{\Gamma}}}
\def\Lap{{\hat{\Lambda}}}
\def\Sip{{\hat{\Sigma}}}
\def\Upsig{{\check{\Upsilon}}}
\def\Kg{{\check{K}}}
\def\ellp{{\hat{\ell}}}
\def\j{j}
\def\jp{{\hat{j}}}
\def\Stir{{\mathscr S}}
\def\BB{{\mathscr B}}
\def\LL{{\mathcal L}}
\def\SS{{\mathscr S}}
\def\DD{{\mathcal D}}
\def\VV{{\mathcal V}}
\def\WW{{\mathcal W}}
\def\OO{{\mathcal O}}
\def\RR{{\mathcal R}}
\def\AA{{\mathscr A}}
\def\CC{{\mathscr C}}
\def\NN{{\mathcal N}}
\def\WW{{\mathcal W}}
\def\HH{{\mathcal H}}
\def\XX{{\mathcal X}}
\def\YY{{\mathcal Y}}
\def\ZZ{{\mathcal Z}}
\def\UU{{\mathcal U}}
\def\XX{{\mathcal X}}
\def\RR{{\mathcal R}}
\def\cir{{\scriptscriptstyle \circ}}
\def\circa{\thickapprox}
\def\vain{\rightarrow}
\def\leqs{\leqslant}
\def\geqs{\geqslant}
\def\ss{s}
\def\vains{\stackrel{\ss}{\rightarrow}}
\def\parn{\par \noindent}
\def\salto{\vskip 0.2truecm \noindent}
\def\spazio{\vskip 0.5truecm \noindent}
\def\vs1{\vskip 1cm \noindent}
\def\fine{\hfill $\diamond$ \vskip 0.2cm \noindent}
\newcommand{\rref}[1]{(\ref{#1})}
\def\beq{\begin{equation}}
\def\feq{\end{equation}}
\def\beqq{\begin{eqnarray}}
\def\feqq{\end{eqnarray}}
\def\barray{\begin{array}}
\def\farray{\end{array}}
%%%%%%%%% THIS NUMBERS EQUATIONS BY SECTIONS %%%%%%%%%%%%%
\makeatletter
\@addtoreset{equation}{section}
\renewcommand{\theequation}{\thesection.\arabic{equation}}
%\thesection instead of \arabic{section} for correct equation numbering
% in appendices
\makeatother
%%%%%%%%%%%%%%%%%%%%%%%%%%%INTESTAZIONE%%%%%%%%%%%%%%%%%%%%%%%%%%%%%%%
\begin{titlepage}
\begin{center}
{\huge On the constants for multiplication in Sobolev spaces.}
\end{center}
\vspace{1truecm}
\begin{center}
{\large
Carlo Morosi${}^1$, Livio Pizzocchero${}^2$} \\
\vspace{0.5truecm}
${}^1$ Dipartimento di Matematica, Politecnico di
Milano, \\ P.za L. da Vinci 32, I-20133 Milano, Italy \\
e--mail: carmor@mate.polimi.it \\
${}^2$ Dipartimento di Matematica, Universit\`a di Milano\\
Via C. Saldini 50, I-20133 Milano, Italy\\
and Istituto Nazionale di Fisica Nucleare, Sezione di Milano, Italy \\
e--mail: livio.pizzocchero@mat.unimi.it
\end{center}
\vspace{1truecm}
\begin{abstract} For $n > d/2$, the Sobolev (Bessel potential) space $H^n(\reali^d, \complessi)$
is known to be a Banach algebra with its standard norm $\|~\|_n$
and the pointwise product; so, there is a best constant $K_{n d}$ such that
$\| f g \|_{n} \leqs K_{n d} \| f \|_{n} \| g \|_{n}$ for all $f, g$
in this space. In this paper we derive upper and lower bounds for these constants,
for any dimension $d$ and any (possibly noninteger) $n \in (d/2, + \infty)$.
Our analysis also includes
the limit cases $n \vain (d/2)^{+}$ and $n \vain + \infty$, for which asymptotic formulas
are presented. Both in these limit cases and for intermediate values of $n$, the lower
bounds are fairly close to the upper bounds. Numerical tables are given for $d=1,2,3,4$,
where the lower bounds are always between $75 \%$ and $88 \%$ of the upper bounds.
\end{abstract}
\vspace{1truecm} \noindent \textbf{Keywords:} Sobolev spaces,
inequalities, pointwise multiplication. \par \vspace{0.4truecm} \noindent \textbf{AMS 2000
Subject classifications:} 46E35, 26D10, 47A60.
\par
\end{titlepage}
\section{Introduction.}
\label{intro}
The theory of Sobolev spaces contains a lot of inequalities
which involve real constants; often, the classical arguments employed to prove these
inequalities allow to infer the existence of such constants, but are unsuitable
to evaluate them accurately. On the other hand, a precise knowledge of
these constants is desirable for several reasons: apart from the intrinsic
interest of the problem, there are many applications where a fully quantitative
analysis relies on these numbers. \parn
The inequality analyzed in this paper refers to the
pointwise multiplication in $H^n(\reali^d, \complessi)$ for any $n > d/2$.
We are interested in the best constant $K_{n d}$ such that
$$ \| f g \|_{n} \leqs K_{n d} \| f \|_{n} \| g \|_{n} $$
for all $f, g \in H^n(\reali^d, \complessi)$, where $\| ~\|_n$ is the standard norm of this space (see
Eq.s \rref{incid} \rref{repfur} later on in this Introduction, and Eq. \rref{bemul} in the next section). \parn
The constants $K_{n d}$ are relevant in relation to
PDEs with polynomial nonlinearities, since they allow
precise estimates on certain approximation methods and on blow up phenomena.
To cite only one example, we refer to the
semilinear heat equation in one space dimension discussed in \cite{approxi};
here, an estimate on $K_{1 1}$ has been employed to
compute the error of the Galerkin approximate solutions,
and the blow up times for certain initial data.
\parn
Evaluating $K_{n d}$ for arbitrary $n$ and $d$ is a nontrivial task. For example, let
the problem be formulated in the variational language: maximize $\| f g \|_n$ with the constraints
$\| f \|_n = \| g \|_n = 1$; if $n$ is integer one can write the corresponding Euler-Lagrange equations,
but these form a cubic system of PDEs of order $2 n$ for $f$ and $g$. \parn
Due to the difficulty of the problem, one could be satisfied even if, in spite of
the exact value of $K_{n d}$, one had sufficiently close lower and upper bounds for it.
Such bounds are proposed in this paper, for any integer $d$ and (possibly noninteger)
$n \in (d/2, +\infty)$. Our upper bounds depend on an accurate use
of the Fourier transform and of the convolution: the conclusion of this analysis is an inequality
$$ K_{n d} \leqs K^{+}_{n d}~, $$
where $K^{+}_{n d}$ is the sup on $[0,+\infty)$ of a function of hypergeometric type. This sup
is easily evaluated, analytically in certain cases and numerically otherwise. \parn
The lower bounds we propose follow directly from the inequality that defines $K_{n d}$,
choosing for $f$ and $g$ appropriate trial functions: these often depend on one or two
real parameters, so one gets the highest lower bound from the chosen functions maximizing
with respect to the parameters. In any case, this procedure gives inequalities of the form
$$ K^{-}_{n d} \leqs K_{n d}~, $$
where $K^{-}_{n d}$ depends on the trial functions: we will consider two specific choices,
giving rise to what we call the "Bessel" or "Fourier" lower bounds. Both types of bounds
are expressible via special functions of hypergeometric type, or by
one-dimensional integrals which are easily computed numerically.
For given values of $n$ and $d$, the best available estimate from below for
$K_{n d}$ is obtained choosing for $K^{-}_{n d}$ the highest between the Bessel and
the Fourier bounds. For
certain values of $n$ only one kind of lower bound is easily computed, so one must be
content with it.
Our investigation
also includes the limit cases $n \vain (d/2)^{+}$ and $n \vain + \infty$; the second limit
requires the asymptotic analysis of certain integrals, which is performed via the Laplace method.
To give an idea of our results, we anticipate some of them. \parn
i) For $n \vain (d/2)^{+}$, it is
$$ K^{+}_{n d} = {\L_d \over \sqrt{n-d/2}} \Big[1 + O(n-{d \over 2}) \Big]~, $$
where $\L_d$ is an explicitly given constant (see the next section, Eq.
\rref{asiep}); on the other hand, denoting with $K^{-}_{n d}$
a conveniently chosen Bessel lower bound, one finds
$$ K^{-}_{n d} = \sqrt{{2 \over 3}}\,
{\L_d \over \sqrt{n-d/2}} \, \Big[1 + O(n-{{d \over 2}}) \Big] ~; $$
so, in this limit $K^{+}_{n d}/K^{-}_{n d} \vain \sqrt{2/3} > 0.816$. \parn
ii) For $n \vain + \infty$, it is
$$ K^{+}_{n d} = T_d
~{ (2/\sqrt{3})^n \over n^{d/4} } \Big[ 1 + O({1 \over n}) \Big]~, $$
with $T_d$ another explicitly given constant (see Eq. \rref{ffas}). On the other hand, denoting with $K^{-}_{n d}$
an appropriate Fourier lower bound, one finds
$$ K^{-}_{n d} = {(5/3)^{1/2} \over 7^{1/4}} T_d
~{ (2/\sqrt{3})^n \over n^{d/4} } \Big[ 1 + O({1 \over n}) \Big]~; $$
thus, $K^{+}_{n d}/K^{-}_{n d} \vain (5/3)^{1/2}  7^{-1/4} > 0.793$. \parn
iii) For $d=1,2,3,4$ we have explored the whole interval $n \in (d/2,+\infty)$, choosing
for each $K^{-}_{n d}$ the most convenient Bessel or Fourier lower bound and comparing it with the upper bound
$K^{+}_{n d}$; for the sample values of $n$ we have considered,
$K^{+}_{n d}/K^{-}_{n d}$ ranges between 0.750 and 0.880. A table of these upper and lower bounds is reported in
the paper. \parn
iv) As previously said, $K^{+}_{n d}$ is the sup of a hypergeometric-like function. Even though this
is easily computed numerically, to avoid this burden one can use
a majorant $K^{++}_{n d} \geqs K^{+}_{n d}$. We define
$K^{++}_{n d}$ using only elementary functions of $n$; this bound
reproduces correctly the asymptotic behavior of $K^{+}_{n d}$ for
$n \vain (d/2)^{+}$, $n \vain + \infty$, and for $1 \leqs d \leqs 7$ is very close to it on the whole range
$(d/2, +\infty)$. \parn
At the end of this Introduction we will give some details on the organization of the paper. Before speaking
about this, we insert a few comments on some related literature.
\vskip 0.2cm \noindent
\textbf{Connections with previous works.} In our paper
\cite{mp2}, we estimated the constants for more general inequalities
related to multiplication in Sobolev spaces; in particular, we discussed
the constants $K_{n a d}$ in the "tame" (or
"Nash-Moser") inequality $$\| f g \|_n \leqs K_{n a d} \max(\| f \|_n \| g \|_a, \| f \|_a \| g \|_n) $$
for $d/2 < a \leqs n$ and $f, g \in H^n(\reali^d, \complessi)$; here $\|~\|_a$ is the norm
of $H^a(\reali^d, \complessi)$. (The cited work is partly related to the previous one \cite{mp1},
and to the subsequent one \cite{funct} on the tame functional calculus in Sobolev spaces).
In the special case $n=a$, the inequality written above coincides with the inequality  of the present paper.
\parn
For arbitrary $d,a,n$, in \cite{mp2} we derived upper and lower bounds
for $K_{n a d}$.
The lower bounds were of the Bessel and Fourier types
also considered here (with no analysis of the limit $n \vain (d/2)^{+}$, and a discussion
of the limit $a$ fixed, $n \vain + \infty$, of course different from the
present limit $n \vain + \infty$; some explicit formulas of \cite{mp2}
for these lower bounds are replaced here with equivalent, but simpler versions,
and we also give some new formula). \parn
The upper bounds for $K_{n a d}$ were obtained by a different method than
the present one for $K_{n d}$; furthermore, if the upper estimates of \cite{mp2} are
applied with $n=a$ they are found to be rougher than the present ones on $K_{n d}$. \parn
The method we use here to get the upper bounds refines an idea which appeared in
\cite{Pos} in relation to the multiplication in the space
$H^n(\bT, \complessi)$, where $\bT := \reali/(2 \pi \interi)$ is
the one-dimensional torus. The author of \cite{Pos} was not
interested in a precise estimate of the constant for multiplication, so he inserted
in his argument some majorization which, although unnecessary,
simplified the proof of the convergence of a series;
the upper bound on the constant for the multiplication in
$H^n(\bT, \complessi)$ arising from this simplification behaves like
const.$\times 2^n$ for large $n$ (see page 294 of the cited paper). Here we replace the one-dimensional torus
with $\reali^d$, and the Fourier series with the $d$-dimensional Fourier transform. The literal
translation of the technique of \cite{Pos} in our framework would give again an upper bound
for $K_{n d}$ behaving like $2^n$ for $n \vain + \infty$; on the contrary, here we use only the strictly necessary
majorizations and finally obtain the bound $K^{+}_{n d}$ involving a hypergeometric function,
which as explained behaves like $(2/\sqrt{3})^n n^{-d/4}$
for $n \vain + \infty$ and is accurate for small $n$ as well.
\vskip 0.2cm \noindent
\textbf{Organization of the paper.} In Section \ref{desc} we state precisely all the results
about the previously mentioned upper and lower bounds for $K_{n d}$.
As a preparation for the proofs, in Section \ref{back} we write a list
of known identities frequently cited in the sequel, on the following subjects: radial integrals, radial Fourier transforms,
hypergeometric functions, integrals with three Bessel functions and the asymptotics of Laplace integrals
(the last two topics are also treated in the Appendices \ref{appint} and \ref{appeaa}).
In Section \ref{upper} we prove all statements
about the upper bounds $K^{+}_{n d}$. In Sections \ref{pbes} and \ref{fourier} we prove
 all the results about the Bessel and Fourier lower bounds, respectively. \parn
In the remaining part of this Introduction, we fix some notations and definitions employed as standards
throughout the paper.
\vskip 0.2cm \noindent
\textbf{Basic notations on $\reali^d$ and Fourier transforms.}
We consider an arbitrary space dimension $d$;
the running variable in $\reali^d$ is $x =
(x_1, ..., x_d)$, and $k = (k_1, ..., k_d)$ when $\reali^d$ is interpreted as
the "wave vector" space of the Fourier
transform. We write $\bullet$ and $|~ |$ for the inner product and the Euclidean norm of $\reali^d$ (so that
$| x | = \sqrt{{x_1}^2 + ... + {x_d}^2}$, $| k | = \sqrt{{k_1}^2 + ... + {k_d}^2}$,
$k \bullet x = k_1 x_1 + ... + k_d x_d$). \parn
We denote with $\FF, \FF^{-1} : S'(\reali^d,
\complessi) \vain S'(\reali^d, \complessi)$ the Fourier
transform of tempered distributions and its inverse, choosing
normalizations so that (for $f$ in $L^1(\reali^d, \complessi)$~)
it is $\FF f(k)= (2 \pi)^{-d/2}$ $\int_{\reali^d} d
x~e^{-i k \bullet x} f(x)$. The restriction of $\FF$ to $\Ld$, with
the standard inner product and the associated norm $\|~\|_{L^2}$,
is a Hilbertian isomorphism. \vskip 0.2cm \noindent
\textbf{Sobolev spaces.} For real $n \geqs 0$, let us introduce the operators
\beq S'(\reali^d, \complessi) \vain S'(\reali^d, \complessi)~,
\qquad g \mapsto \sqrt{1 - \Delta}^{~n}~ g := \FF^{-1} \left(
\sqrt{1 + | \k |^2}^{~n} \FF g \right) \label{lap} \feq
where $\sqrt{1 + | \k |^2}^{~n}$ means the function $k \in \reali^d \mapsto \sqrt{1 + | k |^2}^{~n}$.
The $n$-th order Sobolev (or Bessel potential \cite{Smi}) space of
$L^2$ type and its norm are
\beq \Hn := \lgraffa f \in S'(\reali^d, \complessi)~\Big\vert~
\sqrt{1 - \Delta}^{~n} f \in \Ld~ \rgraffa= \label{incid}\feq
$$ = \lgraffa f \in S'(\reali^d, \complessi)~\Big \vert~
\sqrt{1 + | \k |^2}^{~n}  \FF f \in \Ld \rgraffa~,
$$
\beq \| f \|_{n} := \| \sqrt{1 - \Delta}^{~n}~ f \|_{L^2} =
\|~ \sqrt{1 + | \k |^2}^{~n}~\FF f~ \|_{L^2}~. \label{repfur} \feq
For $n$ integer, these definitions imply
\beq H^n(\reali^d, \complessi) = \{ f \in S'(\reali^d, \complessi)~|~
\nabla^m f \in L^2(\reali^d, \otimes^m \complessi^d)~\forall m
\in \{0, ..., n\}~\}~ \label{hnab} \feq
where
\beq \nabla^m f :=
(\partial_{\lambda_1,..., \lambda_m} f)_{(\lambda_1, ..., \lambda_m) \in \{1, ..., d\}^m}
\label{weput} \feq
and $\partial_{\lambda_i}$ is the distributional derivative with respect to the coordinate
$x_{\lambda_i}$. The statement
$\nabla^m f \in L^2(\reali^d, \otimes^m \complessi^d)$ means that
\beq + \infty > \sum_{\lambda_1,...,\lambda_m =1,...d}
\int_{\reali^d} d x~| (\partial_{\lambda_1,..., \lambda_m} f) (x) |^2
:= \| \nabla^m f \|^{2}_{L^2}~,  \feq
and the norm \rref{repfur} can be written as
\beq \| f \|_n = \sqrt{\sum_{m=0}^n \left( \barray{c} n \\ m \farray \right)
\| \nabla^m f \|^{2}_{L^2}}~. \label{nonab} \feq
\textbf{Other notations. Some useful functions.} The Pochhammer symbol
of $a \in \reali$, $\ell \in \naturali$, is
\beq (a)_\ell := a ( a + 1) ... (a + \ell - 1)~. \label{poch} \feq
The semifactorial of an odd $\m \in \naturali$ is
\beq m!! := 1 . 3 .... (m - 2) m~, \label{semif} \feq
and we also intend $(-1)!! := 1$.
We refer to \cite{Abr} \cite{Luke} \cite{Wat} as our standards for special functions.
In this paper, we frequently use the Gamma function
and its logarithmic derivative $\psi(\z) := \Gamma'(\z)/\Gamma(\z)$;
for future reference, we write here their properties
more frequently employed in the sequel. These are: the shift formulas
\beq \Gamma(\z+1) = \z \Gamma(\z)~, \label{bett} \feq
\beq \psi(\z + 1) = \psi(\z) + {1 \over \z}~; \label{betpsi} \feq
the special values
\beq  \Gamma(1/2) = \sqrt{\pi}~, \Gamma(1) = 1~, \qquad \psi(1/2) = - \ga - 2 \log 2,~
\psi(1) = - \ga~ \label{valgam} \feq
(with $\ga$ the Euler-Mascheroni constant); the duplication formula
\beq \Gamma(2 \z) = {2^{2 \z - 1} \over \sqrt{\pi}}~\Gamma(\z+1/2) \Gamma(\z)~; \label{dupl} \feq
the identity
\beq \int_{0}^{+\infty}
d u~{u^{\sigma-1} \over (1 + u)^{\gamma}} = {\Gamma(\sigma) \Gamma(\gamma-\sigma) \over \Gamma(\gamma)}
\qquad \mbox{for $\gamma > \sigma > 0$}~. \label{recall} \feq
Another function of which we make wide use is
the Gaussian hypergeometric function ${~}_{2} F_{1}(\alpha,\beta,\gamma; \z) \equiv F(\alpha,\beta,\gamma; \z)$.
We are especially interested in the function
\beq F_{n d} : [0,+\infty) \vain (0,+\infty)~, \qquad u \mapsto F_{n d}(u) :=
F\left(2 n - {d \over 2}, n,
n+ {1 \over 2}; - u\right)~, \label{hypf} \feq
$$ d \in \naturali \setminus \{0\}~, \qquad n \in (d/2, + \infty)~. $$
This function has the equivalent representation
\beq F_{n d}(u) = {1 \over (1 + u)^n}~F\left(n,
{d \over 2} + {1 \over 2} - n, n+ {1 \over 2}; {u \over 1 + u} \right)~, \label{hypff} \feq
following from a familiar Kummer transformation (see Sect. \ref{back}, where we return
to some statements appearing here); we also mention the special case
\beq F_{n d}(u) = \!\!\!\!\!\!\sum_{\ell=0}^{n - d/2 - 1/2}
{(n)_{\ell} \, (d/2 + 1/2 - n)_\ell \over (n+1/2)_{\ell} \, \ell!}\, {u^{\ell} \over (1 + u)^{n + \ell}}
\qquad \mbox{for $n - {d \over 2} - {1 \over 2} \in \naturali$}~. \label{semint} \feq
As it often happens dealing with Sobolev spaces,
a central r\^ole in our considerations
is played by the functions
\beq \F_{n d} : \reali^d \vain \complessi~, \qquad k \mapsto \F_{n d}(k) := {1 \over (1 + | k |^2)^n}~; \label{efnd} \feq
\beq \f_{n d} : \reali^d \vain \complessi~, \qquad \f_{n d} := \FF^{-1} \F_{n d}~.  \label{fnd} \feq
It is clear that $\f_{n d} \in \Hn$ if $n > d/2$; explicitly, one has \cite{Smi} \cite{Maz}
\beq \f_{n d}(x) = {| x |^{n - d/2} \over 2^{n - 1} \Gamma(n)}~
K_{n - d/2}(| x |) \label{gemac} \feq
for $x \in \reali^d$; here $K_\nu$ are the modified Bessel functions of the second
kind, or Macdonald functions.
\section{Description of the main results.}
\label{desc}
Let $d \in \naturali \setminus \{0\}$. For (integer or noninteger) $n > d/2$,
the space $\Hn$ is known to be a Banach algebra under the pointwise
multiplication: see, e.g., \cite{Ada}.
\begin{prop}
\label{demul}
\textbf{Definition.} For $n > d/2$, we put
\beq K_{n d} := \min~\{~ K \geqs 0~|~ \| f g \|_n \leqs K \| f \|_n \| g \|_n  ~~
\mbox{for all $f, g \in \Hn$}~\} \label{bemul} \feq
and refer to this as the best (or sharp) constant for the multiplication in $\Hn$.
\end{prop}
In the sequel we present our upper and lower bounds on $K_{n d}$.
\vskip 0.2cm \noindent
\textbf{Upper bounds on} $\boma{K_{n d}}$. These are given by the following proposition, to
be proved in Sect. \ref{upper}.
\begin{prop}
\label{pupper}
\textbf{Proposition.} i) For all $n > d/2$, it is
\beq K_{n d} \leqs K^{+}_{n d} := \sqrt{\sup_{u \in [0,+\infty)} \SS_{n d}(u)}~. \label{kpnd} \feq
\beq \SS_{n d} : [0,+\infty) \vain (0,+\infty)~, \quad
\SS_{n d}(u) := {\Gamma(2 n - d/2) \over (4 \pi)^{d/2} \Gamma(2 n)} (1 + 4 u)^n~F_{n d}(u)~, \label{ff} \feq
with $F_{n d}$ as in Eq. \rref{hypf} or \rref{hypff}.
$\SS_{n d}$ is bounded, and its boundary values for $u = 0$,
$u \vain +\infty$ are
\beq \SS_{n d}(0) = {\Gamma(2 n - d/2) \over (4 \pi)^{d/2} \Gamma(2 n)}~,
\qquad
\SS_{n d}(+\infty) = {\Gamma(n + 1 -d/2) \over 2^{d-1} \pi^{d/2} (n-d/2) \Gamma(n)}~.
\label{limitu} \feq
ii) For $d/2 < n \leqs d/2 + 1/2$ the function $\SS_{n d}$ is increasing, so that
\beq K^{+}_{n d} = \sqrt{\SS_{n d}(+\infty)} =
{1 \over 2^{d/2-1/2} \pi^{d/4}} \sqrt{{\Gamma(n + 1 -d/2) \over (n-d/2) \Gamma(n)} }~.
\label{soth} \feq
For fixed $d$ and $n \vain (d/2)^{+}$, this implies
\beq K^{+}_{n d} = {\L_d \over \sqrt{n-d/2}} \Big[1 + O(n-{d \over 2}) \Big]~,
\qquad \L_d := {1 \over 2^{d/2-1/2} \pi^{d/4} \sqrt{\Gamma(d/2)}}~.
\label{asiep} \feq
iii) For fixed $d$ and $n \vain +\infty$, it is
\beq K^{+}_{n d} = \sqrt{\SS_{n d}({1 \over 2})} \,\Big[1 + O({1 \over n}) \Big] =
T_d \,{ (2/\sqrt{3})^n \over n^{d/4} } \,\Big[1 + O({1 \over n}) \Big],
\quad T_d := {3^{d/4 + 1/4} \over 2^d \pi^{d/4} }.\label{ffas} \feq
\end{prop}
Of course, in Eq.s \rref{limitu} and \rref{soth} we could write
$\Gamma(n+1-d/2)/ (n-d/2)$ $= \Gamma(n-d/2)$;
the expression in the left hand side has been preferred to handle the limit
$n \vain (d/2)^{+}$. Similar choices have been made for other formulas in the sequel.
\vskip 0.2cm \noindent
\textbf{"Bessel" lower bounds on} $\boma{K_{n d}.}$
The general method to obtain lower bounds
on this constant is based on the obvious inequality
\beq K_{n d} \geqs {\| f g \|_n \over \| f \|_n \| g \|_n} \label{ofc} \feq
for all nonzero $f, g \in \Hn$; this gives a lower bound for any pair of "trial functions" $f, g$.
Inspired by \cite{mp2}, we choose for $f$ and $g$ the function
\beq \f_{\la n d}(x) := \f_{n d}(\la x) \label{lagemac} \feq
where $\lambda \in (0, +\infty)$ is a parameter and $\f_{n d}$ is defined by Eq.
\rref{fnd}. By comparison with that equation, we find
\beq \f_{\la n d} = \FF^{-1} \F_{\la n d}~, \qquad \F_{\la n d}(k) :=
{1 \over \lambda^d (1 + | k |^2/\la^2)^n}~. \feq
To give a lower bound for $K_{n d}$ in terms of these
functions simply amounts to compute $\| \f_{\lambda n d} \|_n$, $\|
\f_{\lambda n d}^2 \|_n$. These norms were already calculated in
\cite{mp2}; here we give them in a more simple and complete form, and add
an analysis of the limit case when $n$ is close to $d/2$. All these facts are described
by the forthcoming proposition, to be proved in Sect. \ref{pbes}.
\begin{prop}
\label{pbessel}
\textbf{Proposition.} i) For all $n > d/2$ and $\lambda > 0$, it is
\beq K_{n d} \geqs \Kap^{B}_{n d}(\lambda) := { \| \f^2_{\lambda n d} \|_n \over \| \f_{\lambda n d} \|_n^2}~,
\label{theabove} \feq
whence
\beq K_{n d} \geqs K^{B}_{n d} := \sup_{\lambda > 0} \Kap^{B}_{n d}(\lambda)~. \label{thea1} \feq
The norms in Eq. \rref{theabove} are given by
\beq \| \f_{\la n d} \|^2_{n} = {\pi^{d/2} \Gamma(n + 1 -d/2) \over (n-d/2) \Gamma(n) \, \la^d} ~
F(-n, {d \over 2}, n; 1 - \la^2)~; \label{po} \feq
\beq \| \f_{\la n d} \|^2_{n} = {\pi^{d/2} \over \Gamma(d/2) \Gamma(2 n) \la^d}~\times \label{remar} \feq
$$ \times ~\sum_{\ell=0}^n~ \left( \barray{c} n \\ \ell \farray \right)
\Gamma(\ell + d/2) \Gamma(2 n - d/2 - \ell)~\la^{2 \ell}~
\mbox{for $n$ integer}~; $$
\beq \| \f^2_{\la n d} \|^2_{n} =
{\pi^{d/2} \, \Gamma^2(2 n - d/2) \over \Gamma(d/2) \Gamma^2(2 n) \, \la^d}~
\int_{0}^{+\infty} d u~ u^{d/2-1} (1 + 4 \la^2 u)^n F^2_{n d}(u)~,  \label{poo} \feq
with $F_{n d}$ as in Eq.s (\ref{hypf}-\ref{hypff});
$$ \| \f^2_{\la n d} \|^2_{n} =
{\pi^{d/2} \Gamma^2(2 n - d/2) \over \Gamma(d/2) \Gamma^2(2 n) \la^d}
\sum_{\ell, m=0}^{n - d/2 - 1/2}
{(n)_{\ell} \, (d/2 + 1/2 - n)_\ell   \over (n+1/2)_{\ell} \, \ell!}\,
{(n)_{m} \, (d/2 + 1/2 - n)_m \over (n+1/2)_{m} \, m!} \times $$
$$ \times \, {\Gamma(d/2 + \ell + m) \Gamma(n -d/2) \over \Gamma(n+\ell + m)}~
F(-n, {d \over 2} + \ell + m , n + \ell + m; 1 - 4 \la^2)~$$
\beq \mbox{for $n - {d \over 2} - {1 \over 2}$ integer}~. \label{pooint} \feq
\vskip 0.1cm \noindent
ii) Let $d/2 < n \leqs d/2+1/2$. Then, for all $\lambda > 0$ it is
\beq \| \f^2_{\la n d} \|^2 \geqs \Fs_{n d}(\lambda)~, \label{hgr} \feq
so that
\beq \Kap^{B}_{n d}(\lambda) \geqs
\Kap^{BB}_{n d}(\lambda) := {\sqrt{\Fs_{n d}(\lambda)} \over \| \f_{\lambda n d} \|^2}~, \feq
\beq K^{B}_{n d} \geqs K^{BB}_{n d} := \sup_{\lambda > 0} \Kap^{BB}_{n d}(\lambda)~. \label{kbb} \feq
Here:
\beq \Fs_{n d}(\lambda) :=
{\pi^{d/2} \, \Gamma^2(2 n - d/2) \over (n-d/2)^3 \Gamma^2(2 n) \la^d}\,\Big[
P_{n d}^2 \, {\Gamma(n + 1 -d/2) \over \Gamma(n)} F(-n, {d \over 2}, n; 1 - 4 \lambda^2) +
\label{kapf} \feq
$$ - P_{n d} \, Q_{n d} \, {\Gamma(2 n + 1 - d) \over \Gamma(2 n - d/2)}
\, F(-n, {d \over 2}, 2 n - {d \over 2}; 1 - 4 \lambda^2) + $$
$$ + q^2_{n d} \, {\Gamma(3 n + 1 - 3 d/2) \over
3 \, \Gamma(3 n - d)} \, F(-n, {d \over 2}, 3 n - d; 1 - 4 \lambda^2)~\Big]; $$
\beq P_{n d} := {\Gamma(n+1/2) \Gamma(n + 1 -d/2) \over \sqrt{\pi} \Gamma(2 n - d/2)}~,
\label{pqnd} \feq
$$ Q_{n d} := {\Gamma(n+1/2) \Gamma(d/2 + 1 - n) \over \Gamma(n) \Gamma(1/2 + d/2 - n)}~, \qquad
q_{n d} :=
\left\{ \barray{ll} Q_{n d} & \mbox{if~$P_{n d} \geqs Q_{n d}$,} \\
P_{n d} - (n-d/2) & \mbox{if~ $P_{n d} < Q_{n d}$.} \farray \right. $$
(In the above definition of $Q_{n d}$ one should intend
$\Gamma(0) := \infty$, so that $Q_{n d} = 0$  for $n={d/2} + {1/2}$).
For any fixed $d$, $\lambda$ and for $n \vain (d/2)^{+}$, it is
\beq \Kap^{BB}_{n d}(\lambda) = \sqrt{{2 \over 3}}\,
{\L_d \over \sqrt{n-d/2}} \, \Big[1 + O(n-{{d \over 2}}) \Big]~, \label{sqrt23} \feq
with $\L_d$ as in the asymptotic expression \rref{asiep} for the upper bound $K^{+}_{n d}$
(note that $\sqrt{2/3} > 0.816$).
\end{prop}
\parn
As clarified in the sequel,
the Bessel lower bounds are less interesting for large $n$; therefore, it is not worth to determine
their asymptotics for $n \vain +\infty$.
\vskip 0.2cm \noindent
\textbf{"Fourier" lower bounds on} $\boma{K_{n d}.}$
Another choice for the trial functions amounts to choose for $f$ and $g$ the function
\beq f_{p \sigma d}(x) := e^{i p x_1}~e^{- (\sigma/2) | x |^2} \label{fcar} \feq
where the "Fourier character" $x \vain e^{i p x_1}$ is regularized at
infinity by a Gaussian factor (we take this hint from \cite{mp2}, but we develop it in a
different way).
\parn
As we will see, this choice is especially interesting for large $n$.
The Sobolev norm of any order $n$ of this
function can be expressed using the modified Bessel function of the first kind
$I_{\nu}$, the Pochhammer symbol \rref{poch} and the semifactorial \rref{semif}.
Our results on the Fourier lower bounds
are contained in the forthcoming proposition, to be proved in Sect. \ref{fourier}.
\vskip 0.2cm \noindent
\begin{prop}
\label{pfou} \textbf{Proposition.} i) Let $n > d/2$. For all $p, \sigma >0$, it is
\beq K_{n d} \geqs \Kap^{F}_{n d}(p, \sigma) := { \| f_{2 p, 2 \sigma, d} \|_n \over
\| f_{p \sigma d} \|_n^2}~; \label{thefou} \feq
hence
\beq K_{n d} \geqs K^{F}_{n d} := \sup_{p, \sigma > 0} \Kap^{F}_{n d}(p, \sigma)~. \feq
For all $p, \sigma >0$, it is
\beq \| f_{p \sigma d} \|_n^2 = {2 \, \pi^{d/2}  \over \sigma^{d/2 + 1} p^{d/2 - 1}}~
\int_{0}^{+\infty} d \rho \,\rho^{d/2} (1 + \rho^2)^n e^{-{\rho^2 + p^2\over \sigma}} I_{d/2 - 1} ({2 p \over
\sigma} \rho)~; \label{giveby} \feq
in particular, for $n$ integer it is
$$ \| f_{p \sigma d} \|_n^2 = \pi^{d/2} \sum_{\ell=0}^n \sum_{j=0}^\ell \sum_{g=0}^j
\left( \barray{c} n \\ \ell \farray \right) \left( \barray{c} \ell \\ j \farray \right)
\left( \barray{c} 2 j \\ 2 g \farray \right) {(2 g - 1)!! \over 2^g} \times $$
\beq \times \left({d/2}- {1/2}\right)_{\ell-j}
p^{2 j - 2 g} \sigma^{\ell + g - j - d/2}~. \label{gv} \feq
ii) Fix the attention on the "special" lower bound
\beq K^{FF}_{n d} := \Kap^{F}_{n d}\left(p= {1 \over 2 \sqrt{2}}, \sigma = {3 \over 4 n} \right)~; \label{fix} \feq
then
\beq K^{FF}_{n d} = {(5/3)^{1/2} \over 7^{1/4}} T_d
~{ (2/\sqrt{3})^n \over n^{d/4} } \Big[ 1 + O({1 \over n}) \Big]
\qquad \mbox{for $n \vain + \infty$}~,
\label{ggas} \feq
with $T_d$ as in the asympotic formula \rref{ffas} for the upper bound
(note that $(5/3)^{1/2}/7^{1/4}$ $ > 0.793$).
\end{prop}
\textbf{Remark.} The result \rref{ggas} depends on the
asymptotic analysis of a Laplace integral.
The values for $(p, \sigma)$ in Eq. \rref{fix} have been chosen because they simplify this
analysis, and give rise to the term $(2/\sqrt{3})^n n^{-d/4}$
also appearing in the asymptotics \rref{ffas} for the upper bound.
One could discuss the asymptotics of
$\Kap^{F}_{n d}(p, \sigma =c/n)$ for arbitrary choices
of $p$ and $c$ in $(0,+\infty)$; however, this generalization complicates
the implementation of the Laplace method and, in comparison with \rref{ggas}, yields no sensible increase of
the dominant term.
\vskip 0.2cm \noindent
\textbf{Table of the upper and lower bounds on $\boma{K_{n d}}$ for $\boma{d=1,2,3,4}$ and
some test values of $\boma{n}$.}
This is Table 1, which has been constructed using the upper bounds $K^{+}_{n d}$ given by Prop.
\ref{pupper}, and choosing conveniently one of the lower bounds $K^{B}_{n d}$,
$K^{BB}_{n d}$, $K^{F}_{n d}$, $K^{FF}_{n d}$ in Propositions \ref{pbessel}, \ref{pfou};
the chosen lower bound is generally indicated with $K^{-}_{n d}$, and its type
is specified within the table. We have chosen the values of $n$ within a very wide range,
from $d/2 + 10^{-4}$ to $d/2 + 120$; for a better appreciation of the discrepancy between the
upper and lower bounds, instead of $K^{-}_{n d}$ we have reported the ratio $K^{-}_{n d}/K^{+}_{n d}$. \par \noindent
To compute $K^{+}_{n d}$, we must find the sup of the function $\SS_{n d}$ in Prop. \ref{pupper},
which is given explicitly by item ii) of the same proposition for $d/2 < n \leqs d/2 + 1/2$, and must
be computed directly from the function $\SS_{n d}$ in the other cases; we have done
this numerically in most cases, and sometimes analytically: some examples
are given in Sect. \ref{upper}. For large $n$, the numerical search for the maximum
of $\SS_{n d}$ has been done starting from $u=1/2$, as suggested by
item iii) of Prop \ref{pupper}. \par \noindent
Concerning $K^{-}_{n d}$, we have always chosen for it the most convenient
between the lower bounds in Propositions \ref{pbessel} and \ref{pfou} (i.e.,
the highest one or, in some limit cases, the most easily computable).
\parn
As for the Bessel lower bounds, for $n$ sufficiently distant from $d/2$
we have computed numerically the function $\lambda \vain \Kap^{B}_{n d}(\lambda)$ and its
maximum $K^{B}_{n d}$. For $n$ very close to $d/2$, this computation is
very difficult because the integrals in $\Kap^{B}_{n d}$ converge too slowly;
in this case, we have turned the attention to the function $\lambda \vain \Kap^{BB}_{n d}(\lambda)$ and
estimated numerically its maximum $K^{BB}_{n d}$. \parn
Concerning the Fourier lower bounds, for $n$ not very large we have determined
$K^{F}_{n d}$ maximizing numerically the function $(p, \sigma) \vain \Kap^{F}_{n d}(p, \sigma)$;
for very large $n$, we have turned the attention to the bound $K^{FF}_{n d}$ which is
easily computed numerically. \par \noindent
The Bessel lower bounds are generally higher than the Fourier ones
for small $n$; the contrary happens for large $n$.
\vskip 0.2cm \noindent
\textbf{A more accurate $\boma{n \vain (d/2)^{+}}$ asymptotics for $\boma{K^{+}_{n d}}$.}
This is introduced for the reasons explained in the next paragraph.
For the sake of brevity, let us put
\beq n_d := n - {d \over 2}~; \label{nd} \feq
in place of Eq. \rref{asiep}, we propose a higher order expansion
\beq K^{+}_{n d} = {\L_d \over \sqrt{n_d}} \Big[1 - \M_d \, n_d + O(n_d^2) \Big]~,
\qquad \M_d := {\psi(d/2) + \ga \over 2}~.
\label{asiepp} \feq
This is derived from the explicit expression \rref{soth} of $K^{+}_{n d}$, inserting therein
the expansions
\beq \Gamma(1 + n_d) = \Gamma(1) + \Gamma'(1) \, n_d + O(n_d^2) = 1 - \ga \, n_d + O(n_d^2)~,
\feq
$$ \Gamma(n) = \Gamma({d \over 2}) + \Gamma'({d \over 2}) \, n_d + O(n_d^2) =
\Gamma({d \over 2}) \Big[1 + \psi({d \over 2}) \, n_d + O(n_d^2) \Big]~, $$
(recall that $\Gamma'(\z) = \Gamma(\z) \psi(\z)$, and use Eq. \rref{valgam}).
\vskip 0.1cm \noindent
\textbf{"Elementary" upper bounds $\boma{K^{++}_{n d}}$.}
The results \rref{asiepp} \rref{ffas} on the asymptotics of $K^{+}_{n d}$ in the limits $n\vain (d/2)^{+}$,
$n \vain +\infty$ suggest a way to build new majorants
\beq K^{++}_{n d} \geqs K^{+}_{n d} \geqs K_{n d}~, \qquad n \in (d/2, + \infty)~, \feq
that are presented hereafter.
Even though less precise than the $^{+}$ upper bounds, the $^{++}$ bounds have the
advantage of being elementary functions of $n$; we will show that they are very close to the $^{+}$ bounds
on the whole interval $(d/2, + \infty)$
up to $d = 7$, and fairly close to them up to $d=10$. For any $d$, the elementary $^{++}$ bounds
reproduce the asymptotics \rref{asiepp} \rref{ffas} of the $^{+}$ bounds at the leading order.
\parn
In order to construct $K^{++}_{n d}$, we first
define a function $n \in (d/2, + \infty) \mapsto z_{n d}$ through the equation
$$ K^{+}_{n d} = {(2/\sqrt{3})^n \over n^{d/4}} \Big[ ({3 d \over 8})^{d/4} ~{\L_d \over \sqrt{n_d}}
\left( 1 - {n_d \over n} \right)^{3/2} (1 + \N_d n_d) + T_d ~ \left({n_d \over n}\right)^{3/2} +
z_{n d}~{n_d \over n^2} \Big]~,  $$
\beq \N_d := \log( {\sqrt{3} \over 2} ) + {1 \over 2} + {3 \over d} - \M_d~,
\quad \mbox{$n_d$ as in \rref{nd}}~. \label{eqznd} \feq
This equation is easily solved for $z_{n d}$. From the explicit expression for $z_{n d}$ and from
the asymptotics \rref{asiepp} \rref{ffas}, one gets
\beq z_{n d} = O(\sqrt{n_d}) \quad \mbox{for $n \vain (d/2)^{+}$}~, \qquad z_{n d} = O(1)
\quad \mbox{for $n \vain +\infty$}~; \label{asiznd} \feq
the coefficient $\N_d$ is defined as above just in order to give the first one of these
relations. \parn
On account of Eq.s \rref{asiznd},
for fixed $d$ the function $n \vain z_{n d}$ is bounded on the interval $(d/2, + \infty)$; this ensures
the finiteness of
\beq Z_d := \sup_{n \in (d/2, +\infty)} z_{n d}~. \label{zetd} \feq
Now, putting
\beq
K^{++}_{n d} := {(2/\sqrt{3})^n \over n^{d/4}} \Big[ ({3 d \over 8})^{d/4} ~{\L_d \over \sqrt{n_d}}
\left( 1 - {n_d \over n} \right)^{3/2} (1 + \N_d n_d) + \label{eqkpp} \feq
$$ + T_d ~ \left({n_d \over n}\right)^{3/2} +
Z_{d}~{n_d \over n^2} \Big]~, $$
we see from \rref{eqznd} that $K^{+}_{n d} \leqs K^{++}_{n d}$.
From Eq.s \rref{eqkpp} and
\rref{asiepp} \rref{ffas}, we also infer
\beq {K^{++}_{n d} \over K^{+}_{n d}} = 1 + O(n_d) \quad \mbox{for $n \vain (d/2)^{+}$}~,
\qquad {K^{++}_{n d} \over K^{+}_{n d}} = 1 + O({1 \over n}) \quad \mbox{for $n \vain +\infty$}~. \feq
The forthcoming Table 2 reports, for $1 \leqs d \leqs 10$, the numerical values of the constants $Z_d$ in Eq. \rref{zetd} and
of the quantities
\beq \Theta_d := \sup_{n \in (d/2, +\infty)} {K^{++}_{n d} \over K^{+}_{n d}}~. \feq
The table has been constructed in this way. First of all, for each $d$ in the above range
the function $n \mapsto z_{n d}$ defined by \rref{eqznd} has been plotted (expressing $z_{n d}$
in terms of $K^{+}_{n d}$ and evaluating the latter numerically);
from the graph of $n \mapsto z_{n d}$, the sup $Z_d$ has been evaluated. Secondly,
for the same values of $d$ the ratio $K^{++}_{n d}/K^{+}_{n d}$ has been plotted as a function of $n$,
and its sup $\Theta_d$ has been evaluated from the graph.
\vfill \eject \noindent
{~}
\oddsidemargin=0.3truecm
\vskip -1cm \noindent
\textbf{Table 1. Bounds $\boma{K^{-}_{n d} \leqs K_{n d} \leqs K^{+}_{n d}}$ for $\boma{d=1,2,3,4}$ and}
\textbf{$\boma{n - {d/2} = 10^{-4}}$}, \textbf{$\boma{10^{-2}, 10^{-1}, {1/4}, {1/2}, 1, {3/2}, 3, 6, 15, 30, 60, 120.}$}
\textbf{(The symbol $^{-}$ stands for one
of the types {\STBB, \STB, \STF, \STFF}, indicated below.)}
\vskip 0.3cm \noindent
\hrule
\vskip 0.4cm \noindent
$\boma{d=1}$
\vskip 0.1cm \noindent
\renewcommand\tabcolsep{0.1cm}
\renewcommand\arraystretch{1.3}
{\footnotesize{
$$ \begin{tabular}{c|| c|c|c|c|c|c|c|c|c|c|c|c|c}
$n$ & ${1\over2}\scriptstyle{+10^{-4}}\!$ & ${1 \over 2}\scriptstyle{+10^{-2}}\!$ & ${1 \over 2}\scriptstyle{+10^{-1}}\!$
& $3/4$ & $1$ & $3/2$ & $2$ &
$7/2$ & $13/2$ &
$31/2$ & $61/2$ & $121/2$ & $241/2$ \\[0.1cm]
\hline \hline
$K^{+}_{n d}$ & $56.5$ & 5.69 & 1.90 & 1.30 & 1.00 & 0.852 & 0.814 & 0.834 & 1.07 & 3.09 & 22.4 & 1410 &
$6.63 \! \times \! \! 10^6$ \\ \hline
$\dd{K^{-}_{n d} \over K^{+}_{n d}}$ & $0.816$ & 0.818 & 0.824 & 0.834 & 0.842 & 0.810 & 0.777 & 0.766 & 0.787 & 0.794 & 0.794 & 0.789 &
0.791 \\ [-0.4cm]
$~$ & \TBB & \TBB & \TBB & \TB & \TB & \TB & \TB & \TF & \TF & \TF & \TF & \TFF & \TFF
\end{tabular} $$    }}
\vskip 0.2cm \noindent
$\boma{d=2}$
\vskip 0.1cm \noindent
{\footnotesize{
$$ \begin{tabular}{c|| c|c|c|c|c|c|c|c|c|c|c|c|c}
$n$ & $1\scriptstyle{+10^{-4}}\!$ & $1\scriptstyle{+10^{-2}}\!$ & $1\scriptstyle{+10^{-1}}\!$
& 5/4 & 3/2 & 2 & 5/2 & 4 & 7 & 16 & 31 & 61 & 121
\\[0.1cm] \hline \hline
$K^{+}_{n d}$ & 39.9 & 3.99 & 1.27 & 0.798 & 0.565 & 0.428 & 0.378 & 0.332 & 0.361 & 0.831 & 5.08 & 269 &
$1.07 \! \times \! \! 10^6$ \\ \hline
$\dd{K^{-}_{n d} \over K^{+}_{n d}}$ & 0.816 & 0.817 & 0.826 & 0.844 & 0.865 & 0.842 & 0.811 & 0.752 & 0.772 & 0.788 & 0.794 & 0.786 &
0.789 \\ [-0.4cm]
$~$ & \TBB & \TBB & \TBB & \TB & \TB & \TB & \TB & \TF & \TF & \TF & \TF & \TFF & \TFF
\end{tabular} $$    }}
\vskip 0.2cm \noindent
$\boma{d=3}$
\vskip 0.1cm \noindent
{\footnotesize{
$$ \begin{tabular}{c|| c|c|c|c|c|c|c|c|c|c|c|c|c}
$n$ & ${3 \over 2}\scriptstyle{+10^{-4}}\!$ & ${3 \over 2}\scriptstyle{+10^{-2}}\!$ & ${3 \over 2}\scriptstyle{+10^{-1}}\!$
& 7/4 & 2 & 5/2 & 3 & 9/2 & 15/2 & 33/2 & 63/2 & 123/2 & 243/2
\\[0.1cm] \hline \hline
$K^{+}_{n d}$ & 22.6 & 2.25 & 0.692 & 0.421 & 0.283 & 0.198 & 0.164 & 0.128 & 0.120 & 0.223 & 1.15 & 51.2 &
$1.71 \! \times \! \! 10^5$ \\ \hline
$\dd{K^{-}_{n d} \over K^{+}_{n d}}$ & 0.816 & 0.817 & 0.826 & 0.847 & 0.875 & 0.858 & 0.830 & 0.763 & 0.759 & 0.781 & 0.788 & 0.782 &
0.787 \\ [-0.4cm]
$~$ & \TBB & \TBB & \TBB & \TB & \TB & \TB & \TB & \TB & \TF & \TF & \TF & \TFF & \TFF
\end{tabular} $$    }}
\vskip 0.2cm \noindent
$\boma{d=4}$
\vskip 0.1cm \noindent
{\footnotesize{
$$ \begin{tabular}{c|| c|c|c|c|c|c|c|c|c|c|c|c|c}
$n$ & $2\scriptstyle{+10^{-4}}\!$ & $2\scriptstyle{+ 10^{-2}}\!$ &
$2\scriptstyle{+ 10^{-1}}\!$ & 9/4 & 5/2 & 3 & 7/2 & 5 & 8 & 17 & 32 & 62 & 122
\\[0.1cm] \hline \hline
$K^{+}_{n d}$ & 11.3 & 1.12 & 0.340 & 0.202 & 0.130 & 0.0857 & 0.0678 & 0.0473 & 0.0389 & 0.0590 & 0.259 & 9.72 &
$2.73 \! \times \! \! 10^4$ \\ \hline
$\dd{K^{-}_{n d} \over K^{+}_{n d}}$ & 0.816 & 0.817 & 0.826 & 0.849 & 0.880 & 0.867 & 0.842 & 0.779 & 0.750 & 0.775 & 0.785 & 0.778 &
0.785 \\ [-0.4cm]
$~$ & \TBB & \TBB & \TBB & \TB & \TB & \TB & \TB & \TB & \TF & \TF & \TF & \TFF & \TFF
\end{tabular} $$    }}
\vskip 0.6cm
\noindent
\textbf{Table 2. Constants $\boma{Z_d}$ and $\boma{\Theta_d}$ (for the
elementary upper bounds $\boma{K^{++}_{n d}}$).}
\vskip 0.2cm \noindent
\hrule
\vskip 0.2cm \noindent
{\footnotesize{
$$ \begin{tabular}{c|| c|c|c|c|c|c|c|c|c|c}
$d$ & 1  & 2 & 3 & 4 & 5 & 6 & 7 & 8 & 9 & 10 \\[0.1cm]
\hline \hline
$Z_d$ & 0 & 0.00925 & 0.0458 & 0.0782 & 0.105 & 0.122 & 0.128 & 0.125 & 0.115 & 0.102  \\ \hline
$\Theta_d$ & 1.041 & 1.039 & 1.044 & 1.044 & 1.044 & 1.044 & 1.049 & 1.105 & 1.197 & 1.363  \\ [-0.2cm]
\end{tabular} $$    }}
\vfill \eject \noindent
\oddsidemargin=0.5truecm

\section{Some background.}
\label{back}
In this section we review some known facts, frequently cited in the rest of the paper to prove the
statements of Sect. \ref{desc}.
\vskip 0.2cm \noindent
\textbf{Some $\boma{d}$-dimensional integrals.}
We frequently need to compute integrals of functions on $\reali^d$ which depend only on
the radius $|~|$ (radially symmetric functions), or on the radius and one angle. In this case, we use the formulas
\beq \int_{\reali^d} d x~ \varphi(| x |) = {2 \, \pi^{d/2} \over \Gamma(d/2)}
\int_{0}^{+\infty} d r~ r^{d-1} \varphi(r)~; \label{rag} \feq
\beq \int_{\reali^d} d x~ \chi(| x |, \w \bullet x) = {2 \, \pi^{d/2 - 1/2} \over \Gamma(d/2 - 1/2)}  \times\label{ragan}\feq
$$\times\int_{0}^{+\infty} d r~ r^{d-1} \int_{0}^{\pi} d \theta~ \sin \theta^{d-2} \chi(r, r cos \theta)
\qquad (d \geqs 2; \w \in \reali^d, | \w | = 1)~, $$
holding for all (sufficiently regular) complex valued functions $\varphi$ on $(0,+\infty)$ and $\chi$ on
$(0, +\infty) \times (0, \pi)$. (When writing the analogous formulas for integrals on the
"wave vector" space $(\reali^d, d k)$, the radius $r$ will be renamed $\rho$).
\vskip 0.2cm \noindent
\textbf{Radial Fourier transforms.} Consider two (sufficiently regular) radially symmetric
functions
\beq f : \reali^d \vain \complessi, \quad
x \vain f(x) = \varphi(| x |)~, \qquad F : \reali^d \vain \complessi, \quad
k \vain F(k) = \Phi(| k |)~; \feq
the Fourier and inverse Fourier transforms $\FF f$, $\FF^{-1} F$ are also radially symmetric, and
given by \cite{Boc}
\beq (\FF f)(k) = {1 \over | k |^{d/2 - 1}}~
\int_{0}^{+\infty} d r~ r^{d/2} J_{d/2 - 1}(| k | r) \varphi(r)~, \label{eboc} \feq
\beq (\FF^{-1} F) (x) = {1 \over | x |^{d/2 - 1}}~
\int_{0}^{+\infty} d \rho~ \rho^{d/2} J_{d/2 - 1}(| x | \rho) \Phi(\ro)~, \label{ebocc} \feq
where $J_{\nu}$ are the Bessel functions of the first kind. As anticipated, the latter formula allows
to infer Eq. \rref{gemac} of the Introduction; in this case, Eq.
\rref{ebocc} is applied with $\Phi(\rho) = 1/(1 + \rho^2)^n$ and
the corresponding integral over $\rho$ is given in \cite{Wat}, page 434. \vskip 0.2cm \noindent
\textbf{Hypergeometric function.} As anticipated, in this paper we use extensively the function
$F(\alpha, \beta, \gamma; \x)$; we always interested in real values
of the parameters $\alpha, \beta, \gamma$ and of the argument $\x$.
For future citation, we report here some properties of $F$.
First of all, we cite: the symmetry property
\beq F(\alpha, \beta, \gamma; \x) = F(\beta, \alpha, \gamma; \x)~; \feq
the special values
\beq F(\alpha, \beta, \gamma; 0) = 1~, \qquad \label{specva} \feq
$$ F(\alpha,\beta,\gamma; 1) = {\Gamma(\gamma) \Gamma(\gamma-\alpha-\beta) \over \Gamma(\gamma-\alpha)
\Gamma(\gamma - \beta)}~ \qquad \mbox{for $\gamma > \alpha + \beta$, $\gamma \neq 0,-1,-2,...$}~; $$
the particular cases
\beq F(\alpha, \beta, \beta; \x) = (1-\x)^{-\alpha}~, \label{ey} \feq
\beq F(\alpha, -m, \gamma; \x) = \sum_{\ell=0}^m {(\alpha)_{\ell} (-m)_\ell \over (\gamma)_{\ell}}
\, {\x^{\ell} \over \ell!} \qquad \mbox{for $m \in \naturali$}~. \label{partc} \feq
Secondly, we recall that
\beq F(\alpha, \beta, \gamma; \x) =
{\Gamma(\gamma) \over \Gamma(\beta) \Gamma(\gamma-\beta)}~\int_{0}^{1} d s \, s^{\beta-1} (1-s)^{\gamma-\beta-1}
(1 - \x s)^{-\alpha} > 0 \label{irep} \feq
$$ \quad \mbox{for~ $\gamma > \beta > 0$,~ $\x < 1$}~, $$
\beq F(\alpha, \beta, \gamma; 1-\x) =
{\Gamma(\gamma) \over \Gamma(\beta) \Gamma(\gamma-\beta)}~
\int_{0}^{+\infty} du \, u^{\beta-1} (1 + u)^{\alpha-\gamma} (1 + \x u)^{-\alpha} > 0 \label{gf} \feq
$$ \mbox{for~ $\gamma > \beta > 0$,~ $\x > 0$} $$
(\rref{gf} follows from \rref{irep} with a change of variable $s = u/(1+u)$). \parn
Thirdly, we mention the differentiation formula
\beq {d \over d \x} F(\alpha, \beta, \gamma; \x) = {\alpha \beta \over \gamma}\, F(\alpha+1, \beta+1, \gamma+1; \x)
\label{differ}~; \feq
this formula, combined with the positivity statement in \rref{irep} implies
\beq {d \over d \x} F(\alpha, \beta, \gamma; \x) > 0 \qquad \mbox{for~ $\alpha > 0$,
$\gamma > \beta > 0$,~ $\x < 1$.} \label{incr} \feq
Finally, we report the Kummer transformations
\beq F(\alpha,\beta,\gamma; \x) =
{1 \over (1 - \x)^\beta} F(\beta, \gamma -\alpha, \gamma; {\x \over \x - 1})~, \label{fam} \feq
\beq F(\alpha,\beta, \gamma; \x) = (1-\x)^{\gamma-\alpha-\beta} F(\gamma-\alpha, \gamma-\beta, \gamma; \x)~;
\label{sukum} \feq
the first one allows to pass from the form \rref{hypf} to the form \rref{hypff} for $F_{n d}$. The
positivity of $F_{n d}$ is granted by \rref{irep}.
The expression \rref{semint} of $F_{n d}$ for $n-d/2-1/2$ integer follows from
\rref{hypff} and \rref{partc}.
\vskip 0.2cm \noindent
\textbf{An integral involving Bessel functions.} In Sect. \ref{upper} we will use the integral
\beq I_{\nnu \mmu}(h) := \int_{0}^{+\infty} d r\, r^{\nnu + \mmu + 1} J_{\nnu}(h r) K^2_{\mmu/2}(r)
\qquad (\mu > -1, \nu > 0, h > 0)~, \label{imunu} \feq
involving a Bessel function of the first kind $J_{\mu}$ and the square of a Macdonald function $K_{\nu/2}$. It is
\beq I_{\nnu \mmu}(h) = {\sqrt{\pi} \, \Gamma(\nnu + \mmu + 1) \Gamma(\nnu + \mmu/2 + 1) \over 2^{\nnu + 2}
\Gamma(\nnu + \mmu/2 + 3/2)} \, h^{\nnu} \times \label{result} \feq
$$ \times F(\nnu + \mmu + 1, \nnu + {\mmu/2} + 1, \nnu + {\mmu/2} + {3/2}; - {h^2/4})~. $$
This result is probably known,
but it is not easy to trace it in the most common tables on integrals of Bessel functions; for this reason,
the proof of \rref{result} is given in the Appendix \ref{appint}.
\vskip 0.2cm \noindent
\textbf{Laplace integrals.} The classical theory of these integrals
is widely employed in this paper, to discuss the $n \vain + \infty$ asymptotics of
our bounds on $K_{n d}$. \parn
By a \textsl{standard} Laplace integral, we mean an integral depending on a parameter $n$, of the form
\beq L(n) := \int_{0}^{b} d t~\ppsi(t)~ e^{-n \ffi(t)}~, \label{int} \feq
under the following assumptions:
\beq 0 <n_0 < n < + \infty~, \qquad 0 < b \leqs +\infty~; \label{assgen} \feq
$$  \ffi \in C^1((0, b), \reali)~, \quad \ffi'(t) > 0 \quad \mbox{$\forall~ t \in (0,b)$}~, \quad
\lim_{t \vain 0^{+}} \ffi(t) = 0~, $$
$$  \ppsi \in C((0, b), \reali)~, \qquad
\int_{0}^b d t \, | \ppsi(t) | \, e^{-n \ffi(t)} < + \infty  \quad  \mbox{for all $n$ as above}~.$$
Here and in the sequel, $'$ is the derivative;  we shall also put
\beq \xi := {\ppsi \over \varphi'} \in C((0,b), \reali)~. \feq
The Laplace method gives the $n \vain + \infty$ asymptotics of
$L(n)$, using the idea that the major contributions to
this integral should come from the regions close to the minimum point of $\ffi$, i.e.,
to $t=0$. (In certain cases,
this asymptotics gives a fairly good approximation of $L(n)$
also for non large values of $n$). The asymptotic behavior of $L(n)$ is described by the following proposition
(see, e.g., \cite{Olv}; for uniformity of language, the proof
is reviewed in the Appendix \ref{appeaa}).
\begin{prop}
\label{mainp}
\textbf{Proposition.} Suppose that conditions \rref{assgen} hold, and that
\beq \xi(t) = \sum_{i=0}^{\ell-1} P_i \ffi(t)^{\alpha_i-1} + O(\ffi(t)^{\alpha_{\ell} -1})
\qquad \mbox{for}~~ t \vain 0^{+}, \label{supp} \feq
where $\ell \in \{1,2,...\}$, $P_1,..., P_{\ell-1} \in \reali$,
$0 < \alpha_1 < \alpha_2 < ... < \alpha_{\ell}$. Then
\beq L(n) =
\sum_{i=0}^{\ell-1} P_i {\Gamma(\alpha_i) \over n^{\alpha_{i}}} + O({1 \over n^{\alpha_\ell}})
\qquad \mbox{for $n \vain +\infty$}~.
\label{maineq} \feq
\end{prop}
\textbf{More on Laplace integrals.} By a \textsl{general} Laplace integral, we mean
an integral depending on a parameter $n$ of the form
\beq \Lambda(n) := \int_{a}^{c} d s~\PPsi(s) \, e^{-n \Phi(s)}~, \label{intt} \feq
where
\beq 0 < n_0 < n <  + \infty~, \qquad
- \infty \leqs a < c \leqs +\infty~, \qquad\Phi \in C^1((a, c), \reali) ~, \feq
$$ \PPsi \in C((a, c), \reali), \qquad
\int_{a}^{c} d s~| \PPsi(s) | \, e^{-n \Phi(s)} < + \infty \qquad \mbox{for all $n$ as above}~. $$
Under suitable conditions on $\Phi$, $\Lambda(n)$ can be expressed in terms of one or more
standard Laplace integrals. As a first example, suppose
\beq a > - \infty~, \qquad \Phi'(s) > 0 \quad
\mbox{for all ~ $s \in (a,c)$}~, \qquad \Phi(a) := \lim_{s \vain a^{+}} \Phi(s) > -\infty \feq
(the limit certainly exists by the monotonicity of $\Phi$, but it could be $-\infty$);
then
\beq \Lambda(n) = e^{-n \Phi(a)} L(n)~, \feq
$$ \mbox{$L(n)$ as in \rref{int} with $b := c - a$, ~$\ffi(t) := \Phi(a + t) - \Phi(a)$,~ $\ppsi(t):=\PPsi(a+t)$}~. $$
Similarly, if
\beq c < + \infty~, \qquad \Phi'(s) < 0 \quad
\mbox{for all ~ $s \in (a,c)$}~, \qquad \Phi(c) := \lim_{s \vain c^{-}} \Phi(s) > -\infty~, \feq
we can write
\beq \Lambda(n) = e^{-n \Phi(c)} L(n)~, \label{write} \feq
$$ \mbox{$L(n)$ as in \rref{int} with $b := c - a$, ~$\ffi(t) := \Phi(c - t) - \Phi(c)$, ~ $\ppsi(t):= \PPsi(c-t)$}~. $$
As a final example, suppose
\beq \Phi'(s) \lesseqgtr 0 \quad
\mbox{for~ $s \lesseqgtr h$ \qquad ($h \in (a,c)$)}~; \feq
then we can write
\beq \Lambda(n) = e^{-n \Phi(h)} [L^{-}(n) + L^{+}(n)]~, \label{fitt} \feq
$$ L^{\mp}(n) :=
\int_{0}^{b^{\mp}} d t \, \ppsi^{\mp}(t) e^{- n \ffi^{\mp}(t)}~, \quad
b^{-}:= h - a, ~~~b^{+} := c - h~, $$
$$ \ffi^{\mp}(t) := \Phi(h \mp t) - \Phi(h)~, \quad \ppsi^{\mp}(t) := \PPsi(h \mp t)
\qquad \mbox{for $t \in (0, b^{\mp})$}~,$$
and $L^{\mp}(n)$ are standard Laplace integrals. \parn
In all the previous examples, after reexpressing $\Lambda(n)$ in terms of standard Laplace integrals
one should expand in powers of $\ffi$ or $\ffi^{\mp}$ the functions $\xi := \theta/\ffi'$ or
$\xi^{\mp} := \theta^{\mp}/{\ffi^{\mp}}{\vspace{-0.1cm}'}$. Assuming sufficient smoothness for $\Theta$
and $\Phi$, the coefficients of these expansions can be expressed directly in terms of the
derivatives of $\Theta$ and $\Phi$ at $s=a, c$ or $h$, respectively \cite{Olv}. In the third example, where
$\Phi$ has its minimum at an inner point $h$ of $(a,c)$, there is typically an alternation of equal and
opposite coefficients in the expansions of $\xi^{-}$ and $\xi^{+}$; this yields some cancellation effects
in the expansion of $L^{-}(n) + L^{+}(n)$.
\section{Proofs for the upper bounds on $\boma{K_{n d}.}$}
\label{upper}
Let us write $F \ast G$ for the convolution of two (sufficiently regular) complex functions $F, G$ on $\reali^d$,
given by
\beq (F \ast G)(k) := \int_{\reali^d}
d h~F(k - h) G(h)~.  \feq
We have
\beq \FF (f g) = {1 \over (2 \pi)^{d/2}}~\FF f * \FF g \label{send} \feq
for all sufficiently regular functions
$f$ and $g$ on $\reali^d$ (and in particular, for $f, g$ as in the forthcoming Lemma).
\parn
\begin{prop}
\label{pupper1}
\textbf{Lemma.} For all $n > d/2$, it is
\beq K_{n d} \leqs \sqrt{  \sup_{k \in \reali^d} \ES_{n d}(k) }~, \feq
where
\beq \ES_{n d}(k) := {(1 + | k |^2)^n \over (2 \pi)^d}
\left(\F_{n d} \ast \F_{n d} \right)(k) \label{esnd} \feq
and $\F_{n d}(k) := 1/(1 + | k |^2)^n$ for all $k \in \reali^d$, as in Eq. \rref{efnd}.
\end{prop}
\textbf{Proof.} Consider any two functions $f, g \in \Hn$. Then
\beq \| f g \|_n^2 = \int_{\reali^d} d k (1 + k^2)^n | \FF(f g)(k) |^2 = {1 \over (2 \pi)^{d}}
\int_{\reali^d} d k (1 + k^2)^n | (\FF f \ast  \FF g)(k) |^2~. \label{laprec} \feq
On the other hand, by making explicit the convolution we find
\beq  (\FF f \ast  \FF g)(k)  =  \int_{\reali^d} d h\, \FF f(k - h) \FF g (h)  = \feq
$$ = \int_{\reali^d} d h {1 \over \sqrt{1 + |k - h|^2}^n \sqrt{1 + | h |^2}^n}~
\left( \sqrt{1 + |k - h|^2}^n~\FF f(k - h) \sqrt{1 + |h|^2}^n \FF g (h) \right)~. $$
Now, H\"older's inequality $| \int d h~ U(h) V(h) |^2 \leqs \Big(\int d h | U(h) |^2\Big)
\Big(\int d h~| V(h) |^2 \Big)$ gives
\beq | (\FF f \ast  \FF g)(k) |^2 \leqs C_{n d}(k) P(k)~, \label{dains} \feq
$$ C_{n d}(k) := \int_{\reali^d} {d h \over (1 + | k - h |^2)^n (1 + | h |^2)^n } =
\left(\F_{n d} \ast \F_{n d}\right)(k)~, $$
$$ P(k) := \int_{\reali^d} d h (1 + |k - h|^2)^n | \FF f (k - h) |^2
(1 + | h |^2)^n | \FF g(h) |^2~. $$
Inserting \rref{dains} into Eq. \rref{laprec} we get
\beq \| f g \|_n^2 \leqs {1 \over (2 \pi)^{d}} \int_{\reali^d} d k (1 + | k |^2)^n C_{n d}(k) P(k) \leqs \feq
$$ \leqs \Big(\sup_{k \in \reali^d} {(1 + | k |^2)^n \over (2 \pi)^d} C_{n d}(k) \Big)~\int_{\reali^d} d k P(k)~
= \Big(\sup_{k \in \reali^d} \ES_{n d}(k) \Big) \int_{\reali^d} d k P(k)~. $$
But
\beq \int_{\reali^d} d k \, P(k) = \feq
$$ = \Big(\int_{\reali^d} d k (1 + |k |^2)^n | \FF f(k) |^2 \Big)
\Big(\int_{\reali^d} d h (1 + |h |^2)^n | \FF g(h) |^2 \Big)= \| f \|_n^2 \, \| g \|_n^2~, $$
so we are led to the thesis. \fine
\begin{prop}
\label{pupper2}
\textbf{Lemma.} For $n > d/2$ and $k \in \reali^d$, it is
\beq \ES_{n d}(k) = \SS_{n d}\Big({| k |^2 \over 4}\Big)~, \label{clear} \feq
where $\SS_{n d}$ is the function in Eq. \rref{ff} of Prop. \ref{pupper}.
\end{prop}
\textbf{Proof.} Let us recall that $\F_{n d}$ is the Fourier transform of the function
$\f_{n d}$, already considered in Eq.s \rref{fnd} \rref{gemac}. We have
\beq \ES_{n d}(k) = {(1 + | k |^2)^{n} \over (2 \pi)^d} \left(\FF \f_{n d} \ast \FF \f_{n d} \right)(k)
= {(1 + | k |^2)^{n} \over (2 \pi)^{d/2}} \left(\FF \f^2_{n d} \right)(k)~. \label{es} \feq
But $\f^2_{n d}$ is a radially symmetric function, whose explicit expression in terms of
the Macdonald function is given by \rref{gemac}. We insert this expression in the formula
\rref{eboc} for the radially symmetric Fourier transform and obtain
\beq \left(\FF \f^2_{n d}\right)(k) = {1 \over {2^{2 n - 2} \Gamma^2(n) | k |^{d/2 - 1}}}~
\int_{0}^{+\infty} d r~ r^{2 n - d/2} J_{d/2 - 1}(| k | r)~K^2_{n - d/2}(r)~; \feq
the last integral is computed via Eq. \rref{result}, and the final result is
\beq \left(\FF \f^2_{n d}\right)(k) = {\Gamma(2 n - d/2) \over 2^{d/2} \Gamma(2 n)}
F_{n d}({| k |^2 \over 4})~, \label{rf} \feq
with $F_{n d}$ as in \rref{hypf} or \rref{hypff} (to obtain this, one also uses
Eq. \rref{dupl} for $\Gamma$).
Inserting \rref{rf} into \rref{es} we get the thesis. \fine
\vskip 0.2cm \noindent
\textbf{Proof of Prop \ref{pupper}, item i).} Lemmas \ref{pupper1} and \ref{pupper2} give
immediately the bound \rref{kpnd} for $K_{n d}$, with $\SS$ as in Eq. \rref{ff}. \parn
We now pass to the boundary values of the function $\SS_{n d}$ for $u =0$ and $u \vain + \infty$.
To determine $\SS_{n d}(0)$, use either Eq. \rref{hypf} or Eq. \rref{hypff},
together with Eq. \rref{specva}; the result agrees with Eq. \rref{limitu}.
\parn
To determine $\lim_{u \vain +\infty} \SS_{n d}(u)$ we use Eq. \rref{hypff}, the limits
\beq  \left({1 + 4 u \over 1 + u}\right)^n \vain 4^n~, \qquad  {u \over u + 1} \vain 1 \qquad
\mbox{for $u \vain + \infty$}, \feq
and Eq. \rref{specva}; these relations imply
\beq \SS_{n d}(+\infty) = {2^{2 n - d} \over \pi^{d/2 + 1/2}} {\Gamma(n + 1/2) \Gamma(n-d/2) \over \Gamma(2 n)}
= {\Gamma(n -d/2) \over 2^{d-1} \pi^{d/2} \Gamma(n)}~, \feq
where the last equality follows from \rref{dupl}. This gives
the expression in \rref{limitu} after using \rref{bett} with $\z = n-d/2$. \parn
Of course,
the continuity of $\SS_{n d}$ on $[0,+\infty)$ and the finiteness of its $u
\vain + \infty$ limit ensure that $\SS_{n d}$ is bounded on its domain.
\fine
\textbf{Proof of Prop \ref{pupper}, item ii).}
\textsl{Step 1. The function $\SS_{n d}$ is increasing if $d/2 < n \leqs d/2 + 1/2$}.
To prove this, we use Eq.s \rref{ff} \rref{hypff} and the following remarks: \parn
a) the functions $u \in [0,+\infty) \vain (1+ 4 u)/(1+ u) \in [1,4) $
and $u \in [0,+\infty) \vain u/(1+u) \in [0,1)$ are increasing; \parn
b) the function $\x \in (-\infty,1) \vain F(n, d/2 + 1/2 -n, n+1/2; \x)$ is
increasing for $d/2 < n < d/2 + 1/2$, due to \rref{incr}; in the limit case $n = d/2 +1/2$,
this function equals $1$ everywhere (by  \rref{partc}, with $m=0$). \parn
Of course, the fact that $\SS_{n d}$ is increasing implies $\sup_{[0,+\infty)} \SS_{n d} =
\SS_{n d}(+\infty)$, and this fact, with Eq. \rref{limitu}, yields Eq. \rref{soth}. \parn
\textsl{Step 2. The asymptotics \rref{asiep} of $K^{+}_{n d}$ for $n \vain (d/2)^{+}$}. This
is evident from \rref{soth}.
\fine
Now we must prove item iii) of the same proposition, concerning the $n \vain + \infty$ behavior of $K^{+}_{n d}$;
a fairly long series of Lemmas will be established to this purpose. A main point in this argument
is the integral representation, coming from Eq.s \rref{ff} \rref{hypf} and \rref{irep},
\beq \SS_{n d}(u) = {\Gamma(2 n - d/2) \Gamma(n + 1/2) \over
2^{d} \pi^{d/2 + 1/2} \Gamma(n) \Gamma(2 n)}~ \CC_{n d}(u)~, \label{essg} \feq
$$ \CC_{n d}(u) := (1 + 4 u)^n \int_{0}^{1} d s~{s^{n-1} \over \sqrt{1 - s} \, (1 +  u s)^{2 n - d/2}}~. $$
For future convenience, we write
\beq \CC_{n d}(u) = \AA_{n d}(u) + \BB_{n d}(u)~, \label{gg} \feq
\beq \AA_{n d}(u) := (1 + 4 u)^n \int_{1/4}^{1} d s~{s^{n-1} \over \sqrt{1 - s} \, (1 +  u s)^{2 n - d/2}}~,
\label{aabb} \feq
$$ \BB_{n d}(u) := (1 + 4 u)^n \int_{0}^{1/4} d s~{s^{n-1} \over \sqrt{1 - s} \, (1 + u s)^{2 n - d/2}}~. $$
\parn
\begin{prop}
\label{lemes}
\textbf{Lemma}. Define
\beq B_{n d} := \sup_{u \in [0,+\infty)} \BB_{n d}(u)~; \label{bnd} \feq
then, for fixed $d$ and $n \vain + \infty$,
\beq  B_{n d} = O\Big({1 \over \sqrt{n}} \, ({9 \over 8})^n\Big)~. \label{ogran} \feq
\end{prop}
\textbf{Proof.} We will estimate $\BB_{n d}(u)$ with different methods for $u \in [0,2]$ and
$u \in (2,+\infty)$, respectively. \parn
Let $0 \leqs u \leqs 2$; we reexpress the definition of $\BB_{n d}(u)$ as
\beq \BB_{n d}(u) = (1 + 4 u)^n
\int_{0}^{1/4} d s~{s^{d/4-1} \over \sqrt{1 - s}} \, \left({s \over (1 + u s)^2}\right)^{n - d/4}~. \feq
The function $s \vain s/(1 + u s)^2$ is increasing for $0 \leqs s < 1/u$;
but $1/u > 1/4$, so the maximum of this function for $0 \leqs s \leqs 1/4$ is attained
at $s = 1/4$. From here one gets
\beq \BB_{n d}(u) \leqs (1 + 4 u)^n
\int_{0}^{1/4} d s~{s^{d/4-1} \over \sqrt{1 - s}} \, \left({1/4 \over (1 + u/4)^2}\right)^{n - d/4}~= \feq
$$ = (1 + 4 u)^{d/4} \left((1 + 4 u) \over (2 + u/2)^2 \right)^{n-d/4}
\int_{0}^{1/4} d s~ {s^{d/4 -1} \over \sqrt{1-s}}~. $$
On the other hand, the function $u \vain (1 + 4 u)/(2 + u/2)^2$ is increasing for $0 \leqs u \leqs 2$, and
equals $1$ when $u=2$; from here and from $(1 + 4 u)^{d/4} \leqs 9^{d/4}$ one easily obtains
\beq \sup_{u \in [0,2]} \BB_{n d}(u) \leqs C_d \qquad \mbox{for all $n > d/2$}, \qquad C_d :=
9^{d/4} \int_{0}^{1/4} d s~ {s^{d/4 -1} \over \sqrt{1-s}}~. \label{put1} \feq
We pass to bind $\BB_{n d}$ for $u \in (2, + \infty)$. Returning to Eq. \rref{aabb}, we write
\beq \BB_{n d}(u) \leqs {2 \over \sqrt{3}}
(1 + 4 u)^n \int_{0}^{1/4} d s~{s^{n-1} \over (1 + u s)^{2 n - d/2}} = \feq
$$ = {2 \over \sqrt{3}} {(1 + 4 u)^n \over u^n} \int_{0}^{u/4} d q~{q^{n-1} \over (1 + q)^{2 n - d/2}}~, $$
where the first inequality follows from $1/\sqrt{1-s} \leqs 2/\sqrt{3}$ for $0 \leqs s \leqs 1/4$, and
the subsequent equality is obtained putting $s = q/u$. On the other hand,
$(1 + 4 u)/u < 9/2$ for $u > 2$ and $\int_{0}^{u/4} < \int_{0}^{+\infty}$ on positive
functions, so
\beq \sup_{u \in (2, +\infty)} \BB_{n d}(u) \leqs {2 \over \sqrt{3}} \, ({9 \over 2})^n
\int_{0}^{+\infty} d q {q^{n-1} \over (1 + q)^{2 n - d/2}} = \feq
$$ = {2 \over \sqrt{3}} \, ({9 \over 2})^n {\Gamma(n-d/2) \Gamma(n) \over \Gamma(2 n - d/2)}  $$
(recall Eq. \rref{recall}). We now apply the duplication formula \rref{dupl} with $\z = n-d/4$; this gives
\beq \sup_{u \in (2, +\infty)} \BB_{n d}(u) \leqs
\sqrt{\pi \over 3}~ 2^{2 + d/2}~ ({9 \over 8})^n~ {\Gamma(n-d/2) \over \Gamma(n - d/4)}~
{\Gamma(n) \over \Gamma(n - d/4 + 1/2)} . \label{put2} \feq
Putting together Eq.s \rref{put1} \rref{put2} we get
\beq B_{n d} \leqs
\max \left(C_d, \sqrt{\pi \over 3}~ 2^{2 + d/2}~ ({9 \over 8})^n~ {\Gamma(n-d/2) \over \Gamma(n - d/4)}~
{\Gamma(n) \over \Gamma(n - d/4 + 1/2)}\right) \label{fr} \feq
for all $n > d/2$. As a final step, we recall that (\cite{Olv}, page 119)
\beq {\Gamma(\z + a) \over \Gamma(\z+b)} = \z^{a-b}
\Big[1 + O({1 \over \z}) \Big] \qquad \mbox{for fixed $a, b \in \reali$ and $\z \vain + \infty$}~; \label{ragam} \feq
this implies, for $n \vain + \infty$,
\beq {\Gamma(n-d/2) \over \Gamma(n - d/4)}~
{\Gamma(n) \over \Gamma(n - d/4 + 1/2)} = n^{- 1/2} \Big[1 + O({1 \over n}) \Big] \label{impl} \feq
and Eq.s \rref{fr} \rref{impl} yield the thesis \rref{ogran}. \fine
\begin{prop}
\label{supa}
\textbf{Lemma.} For all $n > d/2$ one has
\beq \sup_{u \in [0,+\infty)} \AA_{n d}(u) \leqs A_{n d}~, \label{bo1} \feq
$$ A_{n d} := 2^{2 n-d/2} {(1 - {d /2 n})^{n-d/2} \over (1 - {d / 4 n})^{2 n-d/2}}
\int_{1/4}^{1} d s \, {1 \over s \sqrt{1-s} \,(4 - s)^{n-d/2}}~. $$
For fixed $d$ and $n \vain + \infty$, it is
\beq A_{n d} = \sqrt{\pi}~ {3^{d/2 + 1/2} \over 2^{d/2} \sqrt{n}}~ ({4 \over 3})^n
\Big[1 + O({1 \over n}) \Big]~. \label{andas} \feq
\end{prop}
\textbf{Proof.} \textsl{Step 1. The bound \rref{bo1}}. The definition
of $\AA_{n d}$ implies
\beq \sup_{u \in [0,+\infty)} \AA_{n d}(u) \leqs \int_{1/4}^{1} d s~
{s^{n-1} \over \sqrt{1-s}} H_{n d}(s)~, \label{thiseq} \feq
$$ H_{n d}(s) := \sup_{u \in [0,+\infty)} {(1 + 4 u)^n \over (1 + s u)^{2 n - d/2}}~. $$
For $s \in (1/4, 1)$, the function
$u \in [0,+\infty) \vain (1 + 4 u)^n /(1 + s u)^{2 n-d/2}$ attains its maximum when $u$ equals
\beq u_{n d}(s) := {8 n + (d - 4 n) s \over 4 (2 n - d) s}~. \feq
Thus
\beq H_{n d}(s) = \left. {(1 + 4 u)^n \over (1 + s u)^{2 n - d/2}} \right|_{u=u_{n d}(s)} =
{(1 - {d \over 2 n})^{n-d/2} \over (1 - {d \over 4 n})^{2 n-d/2}}~ {2^{2 n - d/2} \over  s^n (4 - s)^{n-d/2}}~,
\feq
and inserting this equation into \rref{thiseq} one gets the thesis \rref{bo1}~. \parn
\textsl{Step 2. The asymptotics \rref{andas}.} We reexpress Eq. \rref{bo1} for $A_{n d}$ as
\beq A_{n d} = 2^{2 n-d/2} U_{n d} \int_{1/4}^{1} d s~ \PPsi(s) e^{-(n-d/2) \Phi(s)}~, \label{andre} \feq
$$ U_{n d} := {(1 - {d \over 2 n})^{n-d/2} \over (1 - {d \over 4 n})^{2 n-d/2}}~, \qquad
\PPsi(s) := {1 \over s \sqrt{1-s}}~, \qquad \Phi(s) := \log(4 - s)~. $$
In this representation we recognize a Laplace integral in the parameter $n-d/2$;
it is $\Phi'(s) < 0$ for all $s \in (1/4, 1)$, $\Phi(1) = \log 3$, and
the scheme of Eq.s \rref{intt}-\rref{write}
suggests to rephrase Eq. \rref{andre} as
$$ A_{n d} = 2^{2 n-d/2} U_{n d} \, e^{-(n-d/2) \Phi(1)} L(n-{d \over 2})
= ({3 \over 2})^{d/2}~ ({4 \over 3})^n~ U_{n d} L(n-{d \over 2}) ~, $$
\beq L(m) := \int_{0}^{3/4} d t~\ppsi(t) e^{-m \varphi(t)}~, \label{andlm} \feq
$$ \ppsi(t) := \PPsi(1 - t) = {1 \over \sqrt{t}\, (1-t)}~,
\quad \varphi(t) := \Phi(1-t) - \Phi(1) = \log(1 + {t \over 3})~. $$
The last integral has the standard Laplace form \rref{int}, and the framework of Prop.
\ref{mainp} prescribes to analyze it introducing the function
\beq \xi(t) := {\ppsi(t) \over \varphi'(t)} = {3 + t \over (1 - t) \sqrt{t}}~. \feq
For $t \vain 0^{+}$, one has
$$ \varphi(t) = {t \over 3} + O(t^2)~, \qquad t = 3 \varphi(t) + O(\varphi(t)^2)~, $$
\beq \xi(t) = {3 \over \sqrt{t}} + O(\sqrt{t}) =
{\sqrt{3} \over \sqrt{\varphi(t)}} + O(\sqrt{\varphi(t)})~. \label{concl} \feq
Now, application of Prop. \ref{mainp} to the last relation \rref{concl} gives
\beq L(m) = {\sqrt{3 \pi} \over \sqrt{m}} + O({1 \over m^{3/2}}) \qquad \mbox{for $m \vain + \infty$}~.
\label{lm} \feq
On the other hand (taking the logarithm and expanding)
\beq U_{n d} = 1 + O({1 \over n}) \qquad \mbox{for $n \vain + \infty$}~; \label{und} \feq
inserting Eq.s \rref{lm} \rref{und} into \rref{andlm}, one easily derives the thesis \rref{andas}~. \fine
\begin{prop}
\label{lobg}
\textbf{Lemma.} For fixed $d$ and $n \vain + \infty$, it is
\beq \CC_{n d}({1 \over 2}) = \sqrt{\pi}~ {3^{d/2 + 1/2} \over 2^{d/2} \sqrt{n}}~ ({4 \over 3})^n
\Big[1 + O({1 \over n}) \Big] \label{gndas} \feq
(note that the right hand sides of this equation and \rref{andas} coincide).
\end{prop}
\textbf{Proof.} The definition \rref{essg} gives
\beq \CC_{n d}({1 \over 2}) = 3^n \int_{0}^{1} d s~{s^{n-1} \over \sqrt{1 - s} \, (1 + s/2)^{2 n - d/2}}
= 3^n \int_{0}^{1} d s~\PPsi_d(s) e^{-(n-d/4) \Phi(s)}~, \label{ggnnd} \feq
$$ \PPsi_d(s) := {s^{d/4-1} \over \sqrt{1 - s}}~, \qquad \Phi(s) := 2 \log(1 + s/2) - \log s~. $$
We have again a Laplace integral, with parameter $n - d/4$; one finds
$\Phi'(s) < 0$ for all $s \in (0,1)$, $\Phi(1) = 2 \log(3/2)$ and referring again
to the scheme (\ref{intt}-\ref{write}) we reexpress \rref{ggnnd} as
\beq \CC_{n d}({1 \over 2})= 3^n e^{-(n-d/4) \Phi(1)} L_d(n-{d \over 4})
= ({3 \over 2})^{d/2}~ ({4 \over 3})^n~ L_d(n-{d \over 4}) ~, \label{gdlm} \feq
$$ L_d(m) := \int_{0}^{1} d t~\ppsi_d(t) e^{-m \varphi(t)}~, $$
$$ \ppsi_d(t) := \PPsi_d(1 - t) = {(1 - t)^{d/4-1} \over \sqrt{t}}~,
\quad \varphi(t) := \Phi(1-t) - \Phi(1) = 2 \log(1 - {t \over 3}) - \log(1 - t)~. $$
Following again the scheme of Prop. \ref{mainp}, we introduce the function
\beq \xi_d(t) := {\ppsi_d(t) \over \varphi'(t)} = {(3 - t) (1 - t)^{d/4} \over (1 + t) \sqrt{t}}~. \feq
Let us keep $d$ fixed. It turns out that Eq.s \rref{concl} are again satisfied with the present
choice of $\varphi$ and with $\xi= \xi_d$. Therefore, Prop. \ref{mainp} gives
the asymptotics, analogous to \rref{lm},
\beq L_d(m) = {\sqrt{3 \pi} \over \sqrt{m}} + O({1 \over m^{3/2}}) \qquad \mbox{for $m \vain + \infty$}~;
\label{lmd} \feq
inserting Eq. \rref{lmd} into \rref{gdlm} we obtain the thesis \rref{gndas}. \fine
\begin{prop}
\label{ulem}
\textbf{Lemma.} For fixed $d$ and $n \vain + \infty$, one has
\beq \sup_{u \in [0,+\infty)} \CC_{n d}(u) = \sqrt{\pi}~ {3^{d/2 + 1/2} \over 2^{d/2} \sqrt{n}}~ ({4 \over 3})^n
\Big[1 + O({1 \over n}) \Big] \label{gsupas} \feq
(again, the right hand side is as in Eq. \rref{andas}).
\end{prop}
\textbf{Proof.} We have
\beq \CC_{n d}({1 \over 2}) \leqs \sup_{u \in [0,+\infty)} \CC_{n d}(u) \leqs
A_{n d} + B_{n d}~, \feq
(the upper bound follows from Eq.s \rref{gg} \rref{bnd} and \rref{bo1}).
Both the above bounds on $\sup \CC_{n d}$ have asymptotics as in the right hand side of Eq. \rref{andas}.
For the lower bound, this is granted by Lemma \ref{lobg}. For the upper bound,
this follows from Lemmas \ref{supa} for $A_{n d}$ and \ref{lemes} for $B_{n d}$: the latter
is negligible with respect to the former, since
\beq B_{n d} = A_{n d}~ O\Big( {(9/8)^n \over (4/3)^n} \Big) =
A_{n d}~ O\Big(({27 \over 32})^n\Big) = A_{n d}~ O({1 \over n^\sigma}) \quad \mbox{for any real $\sigma$}~. \feq
\fine
\textbf{Proof of Prop \ref{pupper}, item iii).} Eq. \rref{essg} and
the definition of $K^{+}_{n d}$ in Eq. \rref{kpnd} give
\beq \sqrt{\SS_{n d}({1 \over 2})} = {1 \over 2^{d/2} \pi^{d/4 + 1/4}}
\sqrt{{\Gamma(2 n - d/2) \over \Gamma(2 n)}~{ \Gamma(n + 1/2) \over \Gamma(n)}}~
\sqrt{\CC_{n d}({1 \over 2})}~, \label{givek0} \feq
\beq K^{+}_{n d} = {1 \over 2^{d/2} \pi^{d/4 + 1/4}}
\sqrt{{\Gamma(2 n - d/2) \over \Gamma(2 n)}~{ \Gamma(n + 1/2) \over \Gamma(n)}}~
\sqrt{\sup_{u \in (0, + \infty)} \CC_{n d}(u)}~. \label{givek} \feq
We know that $\CC_{n d}(1/2)$ and $\sup_{[0,+\infty)} \CC_{n d}$
have the same asymptotics up to $O(1/n)$, given by Lemmas \ref{lobg} and \ref{ulem};
furthermore, Eq. \rref{ragam} implies
\beq {\Gamma(2 n - d/2) \over \Gamma(2 n)}= {1 \over (2 n)^{d/2}} \Big[1 + O({1 \over n})\Big]~,
\qquad {\Gamma(n + 1/2) \over \Gamma(n)} = \sqrt{n} \Big[1 + O({1 \over n})\Big] \label{duerap}~, \feq
and inserting these results into Eq.s \rref{givek0} \rref{givek} we obtain the thesis \rref{ffas}.
\fine
\textbf{Computing the upper bounds $\boma{K^{+}_{n d}}\,$.}
a) For $d/2 < n \leqs d/2 + 1/2$, we have for $K^{+}_{n d}$ the explicit
expression \rref{soth}; this was employed to compute the numerical values reported in Table 1 for
these cases. \parn
b) In all the other cases, to compute $K^{+}_{n d}$ one has to maximize the
function $\SS_{n d}$ given by Eq. \rref{ff},
containing the hypergeometric function $F_{n d}$ of Eq.s (\ref{hypf}-\ref{hypff}).
For $n - d/2 - 1/2$ integer, $\SS_{n d}$ has the elementary expression \rref{semint}. \parn
c) Apart from simple exceptions, the maximization of $\SS_{n d}$  must be
performed numerically.
In all the cases analyzed with $n > d/2 +1/2$, we have found numerical evidence (and sometimes an analytical proof)
that $\SS_{n d}$ has a unique maximum point $u = u_{n d} > 1/2$ in the interval $(0,+\infty)$, so that
\beq K^{+}_{n d} = \sqrt{\SS_{n d}(u_{n d})}~. \feq
d) Let us consider, for example, the case $d=2$.
For $n = 2$, Eq.s  \rref{ff} \rref{hypf} give
\beq \SS_{2 2}(u) = {(1 + 4 u)^2 \over 12 \pi} ~F\left(3, 2, {5 \over 2}; -u\right)~; \feq
one finds numerically that $\SS_{2 2}$ attains its maximum at $u_{2 2} \simeq 6.84$.
For $n=5/2$, using \rref{ff} \rref{semint} one finds
\beq \SS_{5/2, 2}(u) = {(1 + 4 u)^{5/2} \over 96 \pi}\, {6 + u \over (1 + u)^{7/2}}~; \feq
the point of absolute maximum of this function is $u_{5/2, 2} = 16/5 = 3.2$, determined analytically by solving
an algebraic equation of second degree. For larger, half-integer values of $n$, $\SS_{n 2}$ is
again elementary, but the analytic determination of its maximum point involves algebraic equations
of order increasing with $n$; thus, a numerical attack is necessary. \parn
Table 1 also considers, for $d=2$, the values $n=4$, $7$, $16$, $31$, $61$, $121$. In all these cases, one finds
numerically a unique maximum point $u_{n 2} \simeq 1.46$,  $0.915$, $0.654$, $0.576$, $0.538$, $0.519$.
Note the approach of this point to the limit value $u=1/2$ for large $n$, as expected from Eq.
\rref{ffas}; due to this behavior, numerical maximization is simple even for very large
values of $n$. \parn
\section{Proofs for the Bessel lower bounds on $\boma{K_{n d}}$.}
\label{pbes}
\textbf{Proof of Prop. \ref{pbessel}, item i)}.
Eq.s (\ref{theabove})-(\ref{thea1}) are obvious; we must
justify the expressions (\ref{po}-\ref{remar}) of $\| \f_{\la n d} \|_n$,
and (\ref{poo}-\ref{pooint}) for $\| \f^2_{\la n d} \|_n$.
\parn
\textsl{Step 1. Computation of $\| \f_{\la n d} \|_n$.} We have
\beq \| \f_{\la n d} \|^2_n = \int_{\reali^d} d k~(1 + | k |^2)^n
| \FF \f_{\la n d} |^2 = {1 \over \lambda^{2 d}}~
\int_{\reali^d} d k~{(1 + | k |^2)^n \over (1 + | k |^2/\la^2)^{2 n}}
\label{thelast} = \feq
$$ = {2 \pi^{d/2} \over \Gamma(d/2) \la^{2 d}}~
\int_{0}^{+\infty} \!\!\!\! d \rho \, \rho^{d-1}~{(1 + \rho^2)^n \over (1 + \rho^2/\la^2)^{2 n}} =
{\pi^{d/2} \over \Gamma(d/2) \la^{d}}~
\int_{0}^{+\infty} \!\!\!\! d u \, u^{d/2-1}~{(1 + \lambda^2 u)^n \over (1 + u)^{2 n}}~. $$
In the last two passages we have used Eq. \rref{rag} for the integral of a radially symmetric function, depending only on
$\rho := | k |$, and then we have changed the
variable to $u = \rho^2/\la^2$. \parn
For $n$ arbitrary, the last integral in $u$ is computed using the identity \rref{gf}; this
gives the thesis \rref{po} (after using Eq. \rref{bett} with $\z = n-d/2$). \parn
For $n$ integer, in the integral over $u$ we
expand $(1 + \la^2 u)^n$ with the binomial formula,
and integrate term by term; this gives Eq. \rref{remar} after treating each term
by \rref{recall}. \parn
\textsl{Step 2. Computation of $\| \f^2_{\la n d} \|_n$.} According to the definition
\rref{lagemac}, the function $\f_{\la n d}$ is obtained
from the $\f_{n d}$ of Eq.\rref{gemac} rescaling by $\lambda$. From here, and
from Eq. \rref{rf} for $\FF \f^2_{n d}$ we infer
\beq \left(\FF \f^2_{\la n d}\right)(k) = {1\over 2^{d/2} \la^{d}}~
{\Gamma(2 n - d/2) \over \Gamma(2 n)}~ F_{n d}({| k |^2  \over 4 \la^2})~,
\label{nonlab} \feq
with $F_{n d}$ as in Eq.s \rref{hypf} or \rref{hypff}; thus,
\beq \| \f^{2}_{\la n d} \|_{n}^2 = \int_{\reali^d} d k~(1 + | k |^2)^n~
| \FF \f^2_{\la n d}(k) |^2 = \feq
$$ = {\pi^{d/2} \Gamma^2(2 n - d/2)\over 2^{d-1} \Gamma(d/2) \Gamma^2(2 n) \,\lambda^{2 d}}~
\int_{0}^{+\infty} \!\!\!\! d\rho \, \rho^{d-1}~(1 + \rho^2)^n
F^2_{n d}({\rho^2 \over 4 \la^2})~. $$
Now, introducing the scaled variable $u := \rho^2/(4\, \la^2) $ we readily obtain the expression \rref{poo} for
$\| \f^2_{\la n d} \|_n$. \parn
Finally, let us consider the case $n - d/2 - 1/2$ integer and show
that Eq. \rref{poo} becomes Eq. \rref{pooint}. In fact, in this case the function
$F_{n d}$ has the elementary expression \rref{semint}; when this is substituted into the integral over $u$ of Eq.
\rref{poo}, we get
\beq \int_{0}^{+\infty} \!\!\!\! d u \,u^{d/2-1} (1 + 4 \la^2 u)^n F^2_{n d}(u) = \feq
$$ = \!\!\!\!\!\!\sum_{\ell, m =0}^{n - d/2 - 1/2}
{(n)_{\ell} \, (d/2 + 1/2 - n)_\ell \over (n+1/2)_{\ell} \, \ell!}\,\,
{(n)_{m} \, (d/2 + 1/2 - n)_\m \over (n+1/2)_{m} \, m!}\, \times $$
$$ \times \, \int_{0}^{+\infty} \!\!\!\! d u \, u^{d/2 +\ell + m - 1} {(1 + 4 \la^2 u)^n  \over (1 + u)^{2 n + \ell + m }}
~; $$
each of the above integrals can be computed  via Eq. \rref{gf}, and the conclusion is the thesis
\rref{pooint}. \fine
To prove the second item in Prop. \ref{pbessel} we need an elementary bound for the
hypergeometric-like function $F_{n d}$, to be substituted in Eq. \rref{poo} for $\| \f^2_{\la n d} \|_n$;
this will require some Lemmas.
\begin{prop}
\textbf{Lemma.}
\label{lemf}
Assume
\beq f \in C([0,1], \reali) \cap C^2([0,1), \reali),\quad R \in C([0,1], \reali) \cap C^1([0,1), \reali),
\quad \ep > 0~; \feq
\beq f'(\x) = (1 - \x)^{\ep-1} \R(\x)~, \quad R'(\x) >0 \qquad \mbox{for $\x \in [0,1)$,} \label{such} \feq
and consider the $C^2$ function
\beq \x \in [0,1) \mapsto  {f(1) - f(\x) \over (1-\x)^{\ep}}~. \label{defr} \feq
Then: \parn
\beq {f(1) - f(\x) \over (1-\x)^{\ep}} \vain {\R(1) \over \epsilon} \qquad \mbox{for $\x \vain 1^{-}$}~, \label{r1} \feq
\beq  {d \over d \x} \, {f(1) - f(\x) \over (1-\x)^{\ep}} > 0 \qquad \mbox{for $\x \in [0,1)$}~. \label{r1m} \feq
The previous facts imply
\beq  f(1) - f(0) <  {f(1) - f(\x) \over (1-\x)^{\ep}} < {\R(1) \over \epsilon} \qquad \mbox{for $\x \in (0,1)$.} \label{r01} \feq
\end{prop}
\textbf{Proof.} By the generalized Lagrange theorem, it is
\beq {F(1) - F(\x) \over G(1) - G(\x)} = {F'(t_\x) \over G'(t_\x)} \quad \mbox{for some $t_\x \in (\x, 1)$}~, \label{lag} \feq
for all $F, G \in C([0,1], \reali) \cap C^1((0,1), \reali)$ with $G'$ never vanishing, and
for all $\x \in [0,1)$. We apply this statement with
\beq f := F~, \qquad G(\x) := - (1 - \x)^{\epsilon}~, \feq
taking into account Eq. \rref{such}; this gives
\beq  {f(1) - f(\x) \over (1-\x)^{\ep}} = {\R(t_\x) \over \epsilon} \qquad \mbox{for $\x \in [0,1)$, with $t_\x \in (\x,1)$,} \feq
and in the limit $\x \vain 1^{-}$ we obtain Eq. \rref{r1}. \parn
In order to prove \rref{r1m}, we observe that
\beq  {d \over d \x} \, {f(1) - f(\x) \over (1-\x)^{\ep}} =
{\epsilon \over (1 - \x)^{\epsilon+1}} \, (f(1) - f(\x)) - {\R(\x) \over 1-\x} \qquad \mbox{for $\x \in [0,1)$.}
\label{defimp} \feq
On the other hand (intending $\int_{\x}^1$ as an improper Riemann integral)
$$ \epsilon \, (f(1) - f(\x)) = \epsilon \int_{\x}^1 d t \, f'(t) =
\epsilon \int_{\x}^1 d t \, (1 - t)^{\epsilon-1} \R(t) = $$
\beq = (1 - \x)^{\epsilon}\, \R(\x)   +
\int_{\x}^1 d t \, (1 - t)^\epsilon R'(t)~. \label{lasteq} \feq
the last equality following from integration by parts. Inserting \rref{lasteq} into \rref{defimp}
we obtain
\beq {d \over d \x} \, {f(1) - f(\x) \over (1-\x)^{\ep}} = {1 \over (1 - \x)^{\epsilon + 1}}
\int_{\x}^1 d t (1 - t)^\epsilon R'(t)~, \feq
and the positivity of $R'$ gives the thesis \rref{r1m}. \parn
Finally the function $\x \in (0,1) \mapsto (f(1) - f(\x))/(1-\x)^\ep$ is increasing, so it is strictly
bounded from below and above by its limits for $\x \vain 0^{+}$ and $\x \vain 1^{-}$; this yields Eq. \rref{r01}. \fine
\begin{prop}
\label{lemh}
\textbf{Lemma.} Let
\beq 0 < a, b < + \infty~; \qquad a+ b < c < a + b + 1~; \qquad \x \in (0,1)~. \label{hypoabc} \feq
Then
\beq 0 < P(a,b,c) - 1 <  {P(a, b, c) - F(a, b, c; \x) \over (1 - \x)^{c-a-b}} <
Q(a,b,c) \label{tesithen} \feq
where
\beq P(a,b,c) := F(a,b,c; 1) = {\Gamma(c) \Gamma(c-a-b) \over \Gamma(c-a) \Gamma(c-b)}~, \feq
$$ Q(a,b,c) := {\Gamma(c) \Gamma(a + b + 1 - c) \over (c-a-b) \Gamma(a) \Gamma(b)}~. $$
\end{prop}
\textbf{Proof.} We apply the previous Lemma with
\beq f := F(a,b,c; .)~, \qquad \epsilon := c - a - b~. \feq
In this case, the differentiation formula \rref{differ} and the subsequent application
of the Kummer transformation \rref{sukum} give
\beq f'(\x)  = (1 - \x)^{\ep-1} \R(\x)~, \qquad \R(\x) := {a b \over c} F(c-a, c-b, c+1; \x)~. \label{abct} \feq
On the other hand, the hypergeometric function
$\x \mapsto F(c-a, c-b, c+1; \x)$ has positive derivative, due to \rref{incr} and to the assumptions
\rref{hypoabc} for $a,b,c$; the same assumptions ensure this function to be continuous also
at $\x=1$, where its value is determined by Eq. \rref{specva}.
Thus all conditions of the previous Lemma are fulfilled
by $f, \epsilon, R$, and Eq. \rref{r01} gives
$$ F(a,b,c; 1) - F(a,b,c; 0) <  {F(a, b, c; 1) - F(a, b, c; \x) \over (1 - \x)^{c-a-b}} < $$
\beq < {a b \over c (c-a-b)} F(c-a,c-b,c+1; 1)~. \label{y1} \feq
But
\beq F(a,b,c; 1) = P(a,b,c), \quad F(a,b,c; 1) - F(a,b,c;0) = P(a,b,c) - 1 > 0~; \feq
the last inequality holds because $F(a,b,c; \cdot)$ is increasing (see
again Eq. \rref{incr}). Finally, the equality
\beq {a b \over c (c-a-b)} F(c-a,c-b,c+1; 1) = Q(a,b,c) \label{y3} \feq
is easily inferred from Eq. \rref{specva}, using the identity \rref{bett}
with $\z = a$ and $\z=c$. Eq.s (\ref{y1})-(\ref{y3}) yield the thesis. \fine
\textbf{Remark.} The idea of employing \rref{sukum} in the above proof has been
suggested by \cite{Pon}, where the usefulness of this transformation has been pointed out in relation
to similar inequalities for $F$. \parn
\begin{prop}
\label{lemk}
\textbf{Lemma}. Let $a, b, c$, $P(a,b,c), Q(a,b,c)$ be as in Lemma \ref{lemh}, and
\beq q(a,b,c) := \left\{ \barray{ll} Q(a,b,c) & \mbox{~~if~~$P(a,b,c) \geqs Q(a,b,c)$,} \\
P(a,b,c) -1 & \mbox{~~if~~ $P(a,b,c) < Q(a,b,c)$.} \farray \right. \feq
Then
\beq F(a,b,c; \x)^2 > P(a,b,c)^2 - 2 P(a,b,c) \, Q(a,b,c) (1-\x)^{c-a-b} + \label{getthe} \feq
$$ + q(a,b,c)^2 \,(1-\x)^{2(c-a-b)} \qquad \mbox{for $\x \in (0,1)$}.  $$
\end{prop}
\textbf{Proof.} \textsl{Step 1. The case $P(a,b,c) \geqs Q(a,b,c)$.} For any $\x \in (0,1)$, the upper bound in Eq.
\rref{tesithen} implies
\beq F(a,b,c; \x) > P(a,b,c) - Q(a,b,c) (1-\x)^{c-a-b}~. \feq
The right hand side in the above equation is positive, so we infer
\beq F(a,b,c; \x)^2 > (P(a,b,c) - Q(a,b,c) (1-\x)^{c-a-b})^2~; \feq
expanding the right hand side we get the thesis \rref{getthe}, since in this case $Q(a,b,c) = q(a,b,c)$. \parn
\textsl{Step 2. The case $P(a,b,c) < Q(a,b,c)$.} We write
\beq F(a,b,c; \x)^2 = [\, P(a,b,c) - (P(a,b,c) - F(a,b,c; \x))\, ]^2 =  \feq
$$ = P(a,b,c)^2 + (P(a,b,c) - F(a,b,c; \x))^2 - 2 P(a,b,c) (P(a,b,c) - F(a,b,c; \x))~. $$
We insert here the bounds on $P(a,b,c) - F(a,b,c; \x)$ coming from
Eq. \rref{tesithen}; this gives
\beq F(a,b,c; \x)^2 > P(a,b,c)^2 + (P(a,b,c)-1)^2 (1 - \x)^{2(c-a-b)} + \feq
$$ - 2 P(a,b,c) Q(a,b,c) (1-\x)^{c-a-b}~, $$
and we have the thesis \rref{getthe} since in this case $q(a,b,c)= P(a,b,c)-1$.
\fine
\textbf{Proof of Prop. \ref{pbessel}, item ii)}. Throughout the proof, $d/2 < n \leqs d/2 + 1/2$. \parn
\textsl{Step 1. For $\x \in (0,1)$ one has}
\beq F\left(n, {d \over 2} + {1 \over 2} - n, n+ {1 \over 2}; \x \right)^2 \geqs \label{firstineq} \feq
$$ \geqs {P_{n d}^2 \over (n-d/2)^2} - {2 P_{n d} Q_{n d} \over (n-d/2)^2}
(1 - \x)^{n-d/2} + {q^2_{n d} \over (n-d/2)^2} (1 - \x)^{2 n-d}~, $$
\textsl{where $P_{n d}, Q_{n d}$ and $q_{n d}$ are as in \rref{pqnd}.}
For $n < d/2 + 1/2$, this follows from application of Lemma \ref{lemk} with $a = n$, $b = d/2 + 1/2-n$,
$c = n + 1/2$; comparing the coefficients in this Lemma with Eq. \rref{pqnd} we see that
\beq P(a,b,c) = {P_{n d} \over n-d/2}~, \qquad  Q(a,b,c) = {Q_{n d} \over n - d/2}~,
\qquad  q(a,b,c) = {q_{n d} \over n-d/2}~. \feq
Let us pass to the limit case $n = d/2 + 1/2$; then, \rref{firstineq} holds as an equality because
$P_{n d} = 1/2$, $Q_{n d}=0$, $q_{n d}=0$, $F\left(n, {d /2} + {1 / 2} - n, n+ {1 / 2}; \x \right) =
F(d/2 + 1/2, 0, d/2+1; \x) = 1$ (by \rref{partc}, with $m=0$). \parn
\textsl{Step 2. Proof of Eq. \rref{hgr}}: $\| \f^2_{\lambda n d} \|^2_n \geqs \Fs_{n d}(\lambda)$,
\textsl{with $\Fs_{n d}(\lambda)$ as in Eq. \rref{kapf}.}
We start from the expression \rref{poo} of $\| \f^2_{\lambda n d} \|^2_n$; the function $F_{n d}$
therein is expressed as in \rref{hypff}, and its square is bounded via
the result of Step 1 (with $\x = u/(1+u)$). This gives
\beq \| \f^2_{\lambda n d} \|^2_n \geqs
{\pi^{d/2} \, \Gamma^2(2 n - d/2) \over (n-d/2)^2 \Gamma(d/2) \Gamma^2(2 n) \la^d}~\times \feq
$$ \times~\int_{0}^{+\infty} du \, u^{d/2-1}
{(1 + 4 \la^2 u)^n \over (1 + u)^{2 n}}
\left( P_{n d}^2 - 2 {P_{n d} \, Q_{n d} \over (1 + u)^{n-d/2}} + {q^2_{n d} \over (1 + u)^{2 n-d}} \right)~.
$$
The above integral can be written as the sum of three integrals of the form \rref{gf};
after computing each of them by \rref{gf}, we apply \rref{bett} with $\z =
n-d/2$, $2 n - d$ and $3 n - 3 d/2$, respectively. The final result is
the minorant for $\| \f^2_{\lambda n d} \|^2_n$ as in Eq.
\rref{kapf}. \parn
\textsl{Step 3. The $n \vain (d/2)^{+}$ limit of $\Kap^{BB}_{n d}(\lambda)$.} Let $d$ and $\lambda \in (0,+\infty)$
be fixed. We start computing the limiting behavior of $\Fs_{n d}(\lambda)$.
For $n \vain (d/2)^{+}$, the coefficients
$P_{n d}, Q_{n d}$ and $q_{n d}$ therein have the same behavior up to $O(n-d/2)$:
\beq P_{n d}, Q_{n d}, q_{n d}  = {\Gamma(d/2+1/2) \over \sqrt{\pi} \Gamma(d/2)}~\Big[1 + O(n-d/2) \Big]~.\label{bev1} \feq
In the same limit, the three hypergeometric functions also have equal behavior:
$$ F(-n, {d \over 2}, n; 1 - 4 \lambda^2),~ F(-n, {d \over 2},  2 n - {d \over 2}; 1 - 4 \lambda^2),~
F(-n, {d \over 2}, 3 n - d; 1 - 4 \lambda^2)~ = $$
\beq = F(-{d \over 2}, {d \over 2}, {d \over 2}; 1 - 4 \lambda^2) + O(n-{d \over 2}) = 2^d \lambda^d +
O(n-{d \over 2})~,
\label{bev2} \feq
where the last equality follows from \rref{ey}.
Inserting Eq.s \rref{bev1} \rref{bev2} into \rref{kapf}, we find
$$ \Fs_{n d}(\lambda) = {2^d \pi^{d/2-1} \, \Gamma(d/2+1/2)^2 \over 3 \, \Gamma^2(d) \Gamma(d/2)}\,
{1 + O(n-d/2) \over (n-d/2)^3}~= $$
\beq = {\pi^{d/2} \over 3\, 2^{d-2} \Gamma(d/2)^3}~{1 + O(n-d/2) \over (n-d/2)^3}~; \label{num} \feq
the second equality in \rref{num} follows from the first one applying the duplication formula
\rref{dupl} with $\z=d/2$. \parn
Let us pass to the $n \vain (d/2)^{+}$ behavior of $\| \f_{\la n d} \|_n$;
from \rref{po} and \rref{ey}, we infer
\beq \| \f_{\la n d} \|^2_n = {\pi^{d/2} \over \Gamma(d/2) \lambda^d} F(-{d \over 2}, {d \over 2}, {d \over 2};
1- \lambda^2)~{1 + O(n-d/2) \over n - d/2} =
\label{den} \feq
$$ = {\pi^{d/2} \over \Gamma(d/2)}~{1 + O(n-d/2) \over n - d/2}~. $$
Since $\Kap^{BB}_{n d}(\lambda)= \sqrt{\Fs_{n d}(\lambda)}/\| \f_{\lambda n d} \|_n^2$, from \rref{num} and
\rref{den} we obtain
\beq  \Kap^{BB}_{n d}(\lambda) =
~{1 \over \sqrt{3} \, 2^{d/2-1} \pi^{d/4} \sqrt{\Gamma(d/2)}} {1 + O(n-d/2) \over \sqrt{n - d/2}}~;
\feq
comparing this with the definition \rref{asiep} of $\L_d$, we get the thesis \rref{sqrt23}. \fine
\textbf{Computing the Bessel lower bounds.} a) For all $n > d/2$, the lower bound $\Kap^{B}_{n d}(\la)$
is the ratio of $\| \f^2_{\la n d} \|_n$ and $\| \f_{\la n d} \|^2_n$. The norm of
$\f_{\la n d}$ has the analytic expression \rref{po} in terms of a hypergeometric
function, that becomes the elementary formula \rref{remar} for $n$ integer. \parn
The norm
of $\f^2_{\la n d}$ has the integral representation \rref{poo}, involving the hypergeometric-like
function $F_{n d}$ of Eq.s (\ref{hypf}) (\ref{hypff}). For $n-d/2-1/2$ integer, this
norm has the explicit expression \rref{pooint} in terms of hypergeometric functions.
For $n - d/2-1/2$ noninteger, the integral in \rref{poo}
must be computed numerically. As anticipated, this is a difficult task for $n$ very close to $d/2$,
due to the slow convergence of the integral:
the integrand behaves like $1/u^{1 + (n - d/2)}$ for $u \vain +\infty$ (as made evident by Eq.
\rref{hypff} for $F_{n d}$), and
we are interested in situations where $n-d/2 = 10^{-4}$. In these cases it is
convenient to compute, in place of $\| \f^2_{\lambda n d} \|^2_n$, the minorant
$\Fs_{n d}(\lambda)$ of Eq. \rref{kapf},
and from this the lower bound $\K^{BB}_{n d}(\lambda)$ of Eq. \rref{kbb}, both of them having analytic
expressions in terms of hypergeometric functions. \parn
b) Assuming we are able to compute
$\Kap^{B}_{n d}(\lambda)$ or $\Kap^{BB}_{n d}(\lambda)$, for each $\lambda$ we have a lower bound
for $K_{n d}$; the next step is maximization with respect to $\lambda$, to get
$K^{B}_{n d}$ or $K^{BB}_{n d}$. In general, this is done
numerically (using some package for automatic maximization or for plotting these functions of
$\lambda$, so as to read the maximum from the graph).
\parn
c) Let us consider, for example, the case $d=2$ and the values of $n$ reported in Table 1.
For $n=3/2$, we have the elementary expression
$$ \Kap^{B}_{3/2, 2}(\lambda) = {\la \over 2 \sqrt{2 \pi}}\,\,
{\sqrt{F(1 - 4 \la^2)} \over F(1 - \la^2)}~,
$$
\beq F(\x) := F(-3/2,1, 3/2; \x) = {5 - 3 \x \over 8} + {3 \over 8} (1 - \x)^2 \EFFE(\x)~, \feq
$$ \EFFE(\x) := \left\{ \barray{lll} \mbox{arctanh}(\sqrt{\x})/\sqrt{\x}~~ & \mbox{if $0 < \x < 1$}~, \\
1 & \mbox{if $\x = 0$}~, \\
\arctan(\sqrt{-\x})/\sqrt{-\x}~~ & \mbox{if $\x < 0$}~. \farray \right. $$
The function $\Kap^{B}_{3/2, 2}$ attains its maximum at $\lambda \simeq 1.38$. $\Kap^{B}_{5/2, 2}$ is also elementary,
with its maximum at $\lambda \simeq 1.36$. For $n = 5/4, 2, 4, 7, 16, 31, 61$
the integral in $\| \f^2_{\lambda n 2} \|_n$ can be computed numerically; from the
graph of $\Kap^{B}_{n 2}$ we have found this function to get its maximum at $\lambda \simeq 1.40, 1.36, 1.39,
1.45, 1.53, 1.57, 1.58$, respectively. \parn
For $n = 1 + 10^{-4}, 1 + 10^{-2}, 1 + 10^{-1}$ the numerical computation of
$\Kap^{B}_{n 2}$ and $K^{B}_{n 2}$ is difficult, so we have turned the attention to the
simpler bound $K^{BB}_{n 2}$; from the analytic expressions of
$\Kap^{BB}_{n 2}(\lambda)$ and numerical optimization, we have found the maximum of this function
to be attained at $\lambda \simeq 1.42$ in each one of the three cases. \parn
For all the cases in the table from $n=5/4$ to $n=61$, the previously mentioned Bessel bounds have been
compared with the Fourier lower bounds $K^{F}_{n 2}$ or $K^{FF}_{n 2}$ of Prop. \ref{pfou}
(for the computation of these Fourier bounds, see the remarks at the end of the following section).
In this way, we have found that the Fourier lower bounds are below the Bessel bounds
up to $n = 5/2$, while the contrary happens for $n > 3/2$ (for example:
$K^{F}_{5/4, 2} < 0.610 \, K^{B}_{5/4, 2}$ and $K^{B}_{61, 2} < 0.411 \, K^{FF}_{61,2}$).
Extrapolating, the Bessel bound $K^{B}_{n 2}$ is likely to be smaller than the
Fourier bounds for the large value $n =121$.
Since the numerical computation of $\Kap^B_{121, 2}$  and $K^{B}_{121, 2}$ is difficult,
in the construction of Table 1 we have chosen directly for $K^{-}_{121, 2}$
a Fourier bound. \parn
\section{Proofs for the Fourier lower bounds on $\boma{K_{n d}.}$}
\label{fourier}
We refer to the trial functions $f_{p \sigma d}$ of Eq. \rref{fcar}. Our aim is to prove all statements
contained in Prop. \ref{pfou}; we will proceed in several steps.
\begin{prop}
\textbf{Lemma.} For all $p, \sigma > 0$ and $n > d/2$, $\| f_{p \sigma d} \|_n$ is given by Eq.
\rref{giveby}.
\end{prop}
\textbf{Proof.} The Fourier transform of $f_{p \sigma d}$ is elementary, and given by
\beq (\FF f_{p \sigma d})(k) = {1 \over \sigma^{d/2}} e^{-{1 \over 2 \sigma} | k - p \w |^2}~,
\qquad \w := (1,0,...,0)~; \feq
thus
\beq \| f_{p \sigma d} \|^2_n = {1 \over \sigma^d} \int_{\reali^d} d k ( 1 + | k |^2)^n
e^{-{1 \over \sigma} | k - p \w |^2} = \label{fff} \feq
$$ = {1 \over \sigma^d} \int_{\reali^d} d k ( 1 + | k |^2)^n
e^{-{| k |^2 + p^2 \over \sigma}  + {2 p \over \sigma} \w \bullet k}~. $$
To go on, let us first consider the case $d=1$. Eq. \rref{fff} gives
\beq \| f_{p \sigma 1} \|^2_n =
{1 \over \sigma} \int_{\reali} d k ( 1 + k^2)^n
e^{-{k^2 + p^2 \over \sigma}  + {2 p \over \sigma} k} =  \feq
$$ = {1 \over \sigma} \int_{0}^{+\infty} d \rho ( 1 + \rho^2)^n
e^{-{\rho^2 + p^2 \over \sigma}} (e^{{2 p \over \sigma} \rho} + e^{- {2 p \over \sigma} \rho})~; $$
(in the last passage, we have used the variable $\rho = | k |$); this gives Eq. \rref{giveby} for $d=1$,
since (\cite{Wat}, page 80)
\beq e^{s} + e^{-s} = \sqrt{2 \pi s} ~ I_{-1/2}(s) \qquad \forall s \in (0,+\infty)~. \feq
Now, let us pass to the case $d \geqs 2$. Eq. \rref{fff}
contains an integral of the form \rref{ragan}, where the integration variable is now $k$ and
$\chi(| k |, \w \bullet k ) = (1 + | k |^2)^n e^{-{1 \over \sigma} (| k |^2+ p^2) + {2 p \over \sigma} \w \bullet k}$; therefore,
\beq \| f_{p \sigma d} \|^2_n =
{2 \pi^{d/2 - 1/2} \over \Gamma({d \over 2} - {1 \over 2}) \sigma^d}
\int_{0}^{+\infty} d \rho~ \rho^{d-1} (1 + \rho^2)^n e^{-{\rho^2 + p^2 \over \sigma}} \int_{0}^{\pi}
d \theta~\sin \theta^{d-2} e^{  {2 p \over \sigma} \rho \cos \theta}~. \feq
On the other hand (\cite{Wat}, page 79)
\beq \int_{0}^{\pi}
d \theta~\sin \theta^{2 \nu} e^{s \cos \theta} = \sqrt{\pi} \,
\Gamma(\nu + {1 \over 2}) \left({2 \over s}\right)^{\nu} I_{\nu}(s)~;
\feq
inserting this result into the previous equation, we obtain the thesis \rref{giveby}. \fine
\begin{prop}
\textbf{Lemma.} For all $p, \sigma > 0$ and integer $n > d/2$, $\| f_{p \sigma d} \|_n$ is given by Eq.
\rref{gv}.
\end{prop}
\textbf{Proof.} We return to the first equation \rref{fff}, and expand $(1 + | k |^2)^n$ by the binomial formula; this gives
\beq \| f_{p \sigma d} \|^2_n =
{1 \over \sigma^d} \sum_{\ell=0}^n \left( \barray{c} n \\ \ell \farray \right)
\int_{\reali^d} d k | k |^{2 \ell} e^{-{1 \over \sigma} | k - p \w |^2}~. \feq
Now, we write the integration variable as $k = (h, q)$,
$(h \in \reali, q \in \reali^{d-1})$; so,
\beq \| f_{p \sigma d} \|^2_n =
{1 \over \sigma^d} \sum_{\ell=0}^n \left( \barray{c} n \\ \ell \farray \right)
\int_{\reali \times \reali^{d-1}} d h\, d q \,(h^2 + | q |^2)^{\ell}\,
e^{-{(h-p)^2 \over \sigma}} \,e^{-{|q |^2 \over \sigma}} =  \label{so} \feq
$$ = {1 \over \sigma^d} \sum_{\ell=0}^n \left( \barray{c} n \\ \ell \farray \right)
\sum_{j=0}^\ell \left( \barray{c} \ell \\ j \farray \right)
\int_{\reali} d h \, h^{2 j} e^{-{(h-p)^2 \over \sigma}}~\int_{\reali^{d-1}} d q | q |^{2 \ell - 2 j}
e^{-{|q |^2 \over \sigma}}~, $$
where, in the last passage, we have used again the binomial formula to expand
$(h^2 + | q |^2)^{\ell}$. On the other hand,
\beq \int_{\reali} d h \,h^{2 j} e^{-{(h-p)^2 \over \sigma}} =
\int_{\reali} d h (h+p)^{2 j} e^{-{h^2 \over \sigma}} = \label{ii1} \feq
$$ = \sum_{m=0}^{2 j} \left( \barray{c} 2 j \\ m \farray \right) p^{2 j - m}
\int_{\reali} d h \,h^{m} e^{-{h^2 \over \sigma}} =
\sum_{g=0}^{j} \left( \barray{c} 2 j \\ 2 g \farray \right) {(2 g - 1)!! \sqrt{\pi}
\over 2^g} p^{2 j - 2 g} \sigma^{1/2 + g}~. $$
The last passage above depends on the evaluation of the integrals with $h^{m}$:
these vanish for $m$ odd, while in the even case $m = 2 g$ we have $\int_{\-\infty}^{+\infty} d h \,
h^{2 g} e^{-h^2/\sigma}$ $= \sigma^{g+1/2} \Gamma(g+ 1/2) $ $= \sigma^{g+1/2} 2^{-g} (2 g - 1)!! \sqrt{\pi}$.
Concerning the integrals over $q$, due to Eq. \rref{rag} we have
$$ \int_{\reali^{d-1}} d q | q |^{2 \ell - 2 j}
e^{-{|q |^2 \over \sigma}} = {2 \pi^{d/2 - 1/2} \over \Gamma(d/2-1/2)}
\int_{0}^{+\infty} d \xi~ \xi^{d - 2 + 2 \ell - 2 j} e^{-{\xi^2 \over \sigma}} = $$
$$ = \pi^{d/2 - 1/2} \sigma^{d/2 - 1/2 + \ell - j}
{\Gamma({d/2} - {1/2} + \ell - j) \over \Gamma(d/2-1/2)} = $$
\beq = \pi^{d/2 - 1/2} \sigma^{d/2 - 1/2 + \ell - j}
(d/2 - 1/2)_{\ell - j}~. \label{ii2} \feq
Inserting Eq.s \rref{ii1} \rref{ii2} into \rref{so}, we finally get the thesis
\rref{gv}. \fine
\textbf{Proof of Prop. \ref{pfou}, item i).} This is given by the two previous Lemmas. \fine
We pass to item ii) of the same proposition, whose proof is more lengthy. The initial step
concerns the expression of $\| f_{p \sigma d} \|_n$ when $p$ is arbitrary
and $\sigma = c/n$ $(c > 0)$; in this case, the already proved Eq. \rref{giveby} becomes
\beq \| f_{p, c/n, d} \|_n^2 = {2 \pi^{d/2} n^{d/2 + 1} \over c^{d/2 + 1} p^{d/2 - 1}}~
\int_{0}^{+\infty} d \rho \, \rho^{d/2} (1 + \rho^2)^n e^{-n {\rho^2 + p^2\over c}} I_{d/2 - 1} ({2 n p \over
c} \rho)~. \label{givvebby} \feq
We will analyze this formula in the limit $n \vain + \infty$. In the first Lemma,
$p$ and $c$ will be arbitrary; in the subsequent ones, based on the theory of
Laplace integrals, we will consider a specific choice,
ultimately yielding Eq. \rref{ggas}.
\parn
\begin{prop}
\textbf{Lemma.} Fix $p > 0$, $c > 0$ and $d$; for $n \vain + \infty$, it is
\beq \| f_{p, c/n, d} \|_n^2 =
{\pi^{d/2 - 1/2} n^{d/2 + 1/2} \over c^{d/2 + 1/2} p^{d/2 - 1/2}} \Big[ \, X_{p c, d/2-1/2}(n) \, + \label{tesbo} \feq
$$ + O\Big({X_{p c, d/2-3/2}(n) \over n}\Big) +  O( (1 + p^2)^n)  \, \Big]~, $$
\beq X_{p c \al}(n) := \int_{p}^{+\infty} d \rho \, \rho^{\al} (1 + \rho^2)^n
e^{-n {(\rho - p)^2 \over c}} \qquad \mbox{for all \, $\al \in \reali$}~. \label{xp} \feq
\end{prop}
\textbf{Proof.} We start from the relations
$$ I_{d/2-1}(s) = {e^{s} \over \sqrt{2 \pi s}} \, h_d(s) = {e^{s} \over \sqrt{2 \pi s}} \, \left(1 + {b_d(s) \over s} \right)
\quad \mbox{for all $s \in (0,+\infty)$}, $$
\beq b_d, h_d  \in L^{\infty}((0,+\infty), \reali)~, \label{wat} \feq
reflecting the asymptotic behavior of the Bessel functions $I_{\nu}(s)$ for $s \vain 0^{+}$
and $s \vain +\infty$ (see \cite{Wat}). \parn
To go on, in Eq. \rref{givvebby} we write $\int_{0}^{+\infty} = \int_{p}^{+\infty} + \int_{0}^p$;
in these two integrals, we substitute the representations \rref{wat} of $I_{d/2-1}$ involving,
respectively, $b_d$ and $h_d$. This gives
\beq \| f_{p, c/n, d} \|_n^2 = {\pi^{d/2 - 1/2} n^{d/2 + 1/2} \over c^{d/2 + 1/2} p^{d/2 - 1/2}}
\Big[ X_{p c, d/2-1/2}(n) + Y_{p c d}(n) + Z_{p c d}(n) \Big]~, \label{intto} \feq
where the $X$ term is defined following Eq. \rref{xp}, and
\beq Y_{p c d}(n) := {c \over 2 p n} \int_{p}^{+\infty} d \rho \, \rho^{d/2-3/2} (1 + \rho^2)^n
e^{-n {(\rho - p)^2 \over c}} b_d({2 p n \over c} \, \rho)~, \feq
$$ Z_{p c d}(n) := \int_{0}^{p}
d \rho \, \rho^{d/2-1/2} (1 + \rho^2)^n e^{-n {(\rho - p)^2 \over c}} h_d({2 n p \over c} \rho)~. $$
We estimate these two integrals. Let $B_d := \sup_{(0,+\infty)} | b_d |$,
$H_d := \sup_{(0,+\infty)} | h_d |$; then
\beq | Y_{p c d}(n) | \leqs {B_d c \over 2 p n} X_{p c, d/2-3/2}(n)~, \feq
$$ | Z_{p c d}(n) | \leqs H_d (1 + p^2)^n \int_{0}^{p} d \rho\, \rho^{d/2-1/2} =
H_d (1 + p^2)^n {p^{d/2+1/2} \over d/2 + 1/2}~, $$
whence
\beq Y_{p c d}(n) = O\Big({X_{p c, d/2-3/2}(n) \over n}\Big)~, \quad Z_{p c d}(n) = O((1 + p^2)^n)
\qquad \mbox{for $n \vain + \infty$}~. \label{sos} \feq
Substituting Eq. \rref{sos} into \rref{intto} we obtain the thesis \rref{tesbo}. \fine
To go on, we observe that Eq. \rref{xp} can be rephrased as
\beq X_{p c \alpha}(n) = \int_{p}^{+\infty} d \rho\, \rho^{\alpha} e^{-n \Phi_{p c}(\rho)}~, \quad
\Phi_{p c}(\rho) := {(\rho - p)^2 \over c} - \log(1 + \rho^2)~. \label{fic}\feq
In the sequel, we apply the Laplace analysis to the integral \rref{fic}. We will consider the
special choice
\beq p := {1 \over 2 \sqrt{2}},~ \quad c = {3 \over 4} \label{spec} \feq
and its double $(2 p, 2 c)$:  this makes easy to compute the minimum point of
$\Phi_{p c}$ and $\Phi_{2 p, 2 c}$. We repeat here the remark made in Sect. \ref{desc},
after stating Prop. \ref{pfou}:
different choices of $(p, c)$ complicate the computations, with no sensible increase
in the dominant term of the Fourier bound $\Kap^{F}(p, c/n)$. (This conclusion is the result of
a tedious analysis, that is not worthy to be reported here).
\begin{prop}
\label{lexp}
\textbf{Lemma.} Let $p, c$ be as in \rref{spec}. For fixed $\alpha \in \reali$ and $n \vain + \infty$, it is
\beq X_{p c \al}(n) =
{3 \sqrt{\pi/5} \over 2^{\al/2 + 1/2} }~{(3/2)^n \over e^{n/6} \sqrt{n}}~
\left[ 1 + O({1 \over n}) \right]~, \label{pc} \feq
\beq X_{2 p, 2 c, \al}(n) = 3 \sqrt{\pi/7}~ 2^{\al/2} {3^n \over e^{n/3} \sqrt{n}}
\left[ 1 + O({1 \over n}) \right] . \label{2pc} \feq
\end{prop}
\textbf{Proof.} \textsl{Step 1. Proof of Eq. \rref{pc}.} We  put for brevity
\beq X_{\alpha}(n) := X_{p c \alpha}(n)~, \qquad \Phi := \Phi_{p c}~. \feq
Explicitly
\beq \Phi(\rho) = {4 \over 3} (\rho - {1 \over 2 \sqrt{2}})^2 - \log(1 + \rho^2)~; \feq
it is easily checked that
$$ \Phi'(\rho) = {2 \over 3} (\rho - {1 \over \sqrt{2}}) {4 \rho^2 + \sqrt{2} \rho + 2 \over 1 + \rho^2}
\lesseqqgtr 0~~\mbox{for}~~\rho \lesseqqgtr {1 \over \sqrt{2}}~, $$
\beq \Phi({1 \over \sqrt{2}}) = {1 \over 6} - \log({3 \over 2})~. \label{check} \feq
Now, following the scheme of \rref{fitt} we reexpress the integral under examination as
\beq X_{\al}(n) = e^{- n \Phi(1/\sqrt{2})} [ L^{-}_{\al}(n) + L^{+}_{\al}(n) ] = {(3/2)^n \over e^{n/6}}
~[ L^{-}_{\al}(n) + L^{+}_{\al}(n) ]~, \label{subint} \feq
$$ L^{-}_{\al}(n) := \int_{0}^{1/(2 \sqrt{2})} d t~\ppsi^{-}_{\al}(t) e^{-n \ffi^{-}(t)}~,
\quad L^{+}_{\al}(n) := \int_{0}^{+\infty} d t~\ppsi^{+}_{\al}(t) e^{-n \ffi^{+}(t)}~; $$
\beq \ppsi^{\mp}_{\al}(t) := ({1 \over \sqrt{2}} \mp t)^{\al}~, \feq
$$ \ffi^{\mp}(t) := \Phi({1 \over \sqrt{2}} \mp t) -
\Phi({1 \over \sqrt{2}}) = \mp {2 \sqrt{2} \over 3} t +
{4 \over 3} t^2 - \log(1 \mp {2 \sqrt{2} \over 3} t + {2 \over 3} t^2)~. $$
The above two integrals have the standard Laplace form discussed in Prop. \ref{mainp}.
Following the usual scheme, we fix the attention on the functions
\beq \xi^{\mp}_{\al}(t) := {\ppsi^{\mp}_{\al}(t) \over {\ffi^{\mp}}{\vspace{-0.1cm}'}(t)} =
{3 \over 4 t}~ {3 \mp 2 \sqrt{2} t + 2 t^2 \over 5 \mp 5 \sqrt{2} t + 4 t^2} ({1 \over \sqrt{2}} \mp t)^{\al}
~. \feq
For $t \vain 0^{+}$, one easily checks that
\beq \ffi^{\mp}(t) = {10  \over 9} \, t^2
\mp { 20 \sqrt{2} \over 81} \, t^3 + O(t^4)~, \label{viaeq} \feq
$$ t = {3 \over \sqrt{10}} \sqrt{\ffi^{\mp}(t)}
\pm {\sqrt{2} \over 10} \ffi^{\mp}(t) + O(\ffi^{\mp}(t)^{3/2})~; $$
\beq \xi^{\mp}_{\al}(t) = {1 \over 2^{\al/2}}\Big[{9 \over 20 t}
\mp  {1 \over \sqrt{2}} ({9 \al \over 10} - {3 \over 10})~\Big] + O(t) = \label{identit} \feq
$$ = {1 \over 2^{\al/2 + 1/2}} \Big[ \, {3 \over 2 \sqrt{5} \, \sqrt{\ffi^{\mp}(t)}} \mp
({9 \al \over 10} - {1 \over 5}) \Big] + O(\sqrt{\ffi^{\mp}(t)})~. $$
We can now apply Prop. \ref{mainp} to both integrals $L^{\mp}_{\al}(n)$; this gives
\beq L^{\mp}_{\al}(n) = {1 \over 2^{\al/2 +1/2}} \Big[ \, {3 \sqrt{\pi} \over 2 \sqrt{5} \, \sqrt{n}}
\mp ({9 \al \over 10} - {1 \over 5}) \, {1 \over n}
\, \Big] + O({1 \over n^{3/2}}) \quad \mbox{for $n \vain + \infty$}~, \feq
and substituting these expansions into Eq. \rref{subint} we get the thesis \rref{pc}. (Note
the mutual cancellation of the terms $\mp (9 \alpha/10 - 1/5)(1/n)$,
in agreement with the remark concluding Sect. \ref{back}). \parn
\textsl{Step 2. Proof of Eq. \rref{2pc}}. In this case, we put
\beq X_{\alpha}(n) := X_{2 p, 2 c, \alpha}(n)~, \qquad \Phi := \Phi_{2 p, 2 c}~. \feq
One has
\beq \Phi(\rho) = {2 \over 3} (\rho - {1 \over \sqrt{2}})^2 - \log(1 + \rho^2)~; \feq
$$ \Phi'(\rho) = {2 \over 3} (\rho - \sqrt{2}) {2 \rho^2 + \sqrt{2} \rho + 1 \over 1 + \rho^2}
\lesseqqgtr 0~~\mbox{for}~~\rho \lesseqqgtr \sqrt{2}~, $$
$$ \Phi(\sqrt{2}) = {1 \over 3} - \log 3~. $$
We can write
\beq X_{\al}(n) = e^{- n \Phi(\sqrt{2})} [ L^{-}_{\al}(n) + L^{+}_{\al}(n) ] = {3^n \over e^{n/3}}
~[ L^{-}_{\al}(n) + L^{+}_{\al}(n) ]~, \label{subint2} \feq
$$ L^{-}_{\al}(n) := \int_{0}^{1/\sqrt{2}} d t~\ppsi^{-}_{\al}(t) e^{-n \ffi^{-}(t)}~,
\quad L^{+}_{\al}(n) := \int_{0}^{+\infty} d t~\ppsi^{+}_{\al}(t) e^{-n \ffi^{+}(t)}~; $$
\beq \ppsi^{\mp}_{\al}(t) := ({\sqrt{2}} \mp t)^{\al}~, \feq
$$ \ffi^{\mp}(t) := \Phi(\sqrt{2} \mp t) - \Phi(\sqrt{2}) = \mp {2 \sqrt{2} \over 3} t +
{2 \over 3} t^2 - \log(1 \mp {2 \sqrt{2} \over 3} t + {1 \over 3} t^2)~. $$
We introduce the functions
\beq \xi^{\mp}_{\al}(t) := {\ppsi^{\mp}_{\al}(t) \over {\ffi^{\mp}}{\vspace{-0.1cm}'}(t)} =
{3 \over 2 t}~ {3 \mp 2 \sqrt{2} t + t^2 \over 7 \mp 5 \sqrt{2} t + 2 t^2} (\sqrt{2} \mp t)^{\al}
~. \feq
For $t \vain 0^{+}$, comparing the expansions of $\ffi^{\mp}$, $\xi^{\mp}_{\al}$ in powers of $t$
we get
\beq \xi^{\mp}_{\al}(t) = 2^{\al/2} \Big[ \, {3 \over 2 \sqrt{7} \, \sqrt{\ffi^{\mp}(t)}}
\mp {1 \over \sqrt{2}} \, ( {9 \al \over 14} - {2 \over 49}) \, \Big] + O(\sqrt{\ffi^{\mp}(t)})~.
\label{identit2} \feq
Applying Prop. \ref{mainp} to $L^{\mp}_{\al}(n)$ we obtain
\beq L^{\mp}_{\al}(n) = 2^{\al/2} \Big[ \, {3 \sqrt{\pi} \over 2 \sqrt{7} \, \sqrt{n}}
\mp {1 \over \sqrt{2}} \, ( {9 \al \over 14} - {2 \over 49}) \, {1 \over n} \, \Big]
+ O({1 \over n^{3/2}}) \quad \mbox{for $n \vain + \infty$}~, \feq
and substituting these expansions into Eq. \rref{subint2} we get the thesis \rref{2pc}. \fine
\vskip 0.2cm \noindent
\begin{prop}
\label{lemr}
\textbf{Lemma.} Let $p, c$ be as in \rref{spec}.  For fixed $d$ and $n \vain + \infty$, it is
\beq \| f_{p, c/n, d} \|^2_n = {2^{3 d/2} \pi^{d/2} \over 3^{d/2 - 1/2} \sqrt{5}}~
{(3/2)^n \over e^{n/6}} \, n^{d/2} \, \left[ 1 + O({1 \over n}) \right] , \label{lapri} \feq
\beq \| f_{2 p, 2 c/n, d} \|^2_n =
{2^d \pi^{d/2} \over 3^{d/2 - 1/2} \sqrt{7}}~
{3^n \over e^{n/3}} \, n^{d/2} \, \left[ 1 + O({1 \over n}) \right]. \label{lasec} \feq
\end{prop}
\textbf{Proof.} To prove Eq. \rref{lapri}, we note that \rref{pc} implies
\beq X_{p c, d/2-1/2}(n) =
{3 \sqrt{\pi/5} \over 2^{d/4 + 1/4} }~{(3/2)^n \over e^{n/6} \sqrt{n}}~
\left[ 1 + O({1 \over n}) \right]~, \label{pcd} \feq
$$ {X_{p c, d/2-3/2}(n) \over n} = {(3/2)^n \over e^{n/6} \sqrt{n}}~
O({1 \over n})~. $$
We insert these results into Eq. \rref{tesbo} for $\| f_{p, c/n, d} \|^2_n$, taking into account
that the present choices of $p, c$ imply
$$ c^{d/2 + 1/2} p^{d/2-1/2} = {3^{d/2 + 1/2} \over 2^{7 d/4 + 1/4}}~; \quad
1 + p^2 = {9 \over 8} = {3/2 \over e^{1/6}} \, \theta, \quad 0.8 <\theta < 0.9 ~; $$
\beq (1 + p^2)^n = {(3/2)^n \over e^{n/6}}~\theta^n = {(3/2)^n \over e^{n/6} \sqrt{n}}~O({1 \over n})~. \feq
The proof of Eq. \rref{lasec} is very similar, depending on Eq.s \rref{2pc} \rref{tesbo}~. \fine
\par \noindent
{\textbf{Proof of Prop. \ref{pfou}, item ii).} This item concerns the $n \vain +\infty$ limit
for the special Fourier lower bound $K^{FF}_{n d}$; comparing the definition \rref{fix} of this bound
with the notations of this section, we see that
\beq K^{FF}_{n d} = {\| f_{2 p, 2 c/n, d} \|_n \over \| f_{p, c/n, d} \|_n^2}~, \qquad
\mbox{$(p, c)$ as in \rref{spec}}~. \feq
From Eq.s \rref{lapri} \rref{lasec} we infer, for $n \vain + \infty$,
\beq K^{FF}_{n d} =
{\sqrt{5} \over 7^{1/4}}~{3^{d/4 - 1/4} \over 2^d \pi^{d/4}}~{(2/\sqrt{3})^n \over n^{d/4}}~
{\dd \sqrt{ 1 + O({1 \over n}) }  \over \dd 1 + O({1 \over n}) } = \feq
$$ = {(5/3)^{1/2} \over 7^{1/4}} T_{d} ~{(2/\sqrt{3})^n \over n^{d/4}}~\Big[ 1 + O({1 \over n})~\Big]. $$
In the last passage we have used the definition \rref{ffas} of $T_d$;
our result is just the thesis \rref{ggas}. \fine
\textbf{Computing the Fourier lower bounds.} a) For any $n$ and $d$, the function
$(p, \sigma) \vain \Kap^{F}_{n d}(p, \sigma)$ in Eq. \rref{thefou} is determined by the function
$(p, \sigma) \vain \| f_{p \sigma d} \|_n$. For $n$ noninteger and given $(p, \sigma)$,
this can be computed
via Eq. \rref{giveby}, evaluating numerically the integral therein; for $n$
integer, we have the elementary expression \rref{gv}. \parn
The bound $K^{F}_{n d}$ is obtained
maximizing $\Kap^{F}_{n d}(p, \sigma)$ with respect to
$(p, \sigma) \in (0,+\infty)^2$; in typical situations this must be done numerically, even for integer $n$
(in any case, the maximization problem is not dramatic because $\Kap^{F}_{n d}(p, \sigma)$ is a lower
bound for \textsl{all} choices of $(p, \sigma)$, even not close to the maximizing pair). \parn
For very large values of $n$, instead of maximizing $\Kap^{F}_{n d}(p, \sigma)$ one can
evaluate it at $(p, \sigma) = (1/(2 \sqrt{2}), 3/(4 n))$; this yields the
the bound $K^{FF}_{n d}$ of Eq. \rref{fix}, that we know to be effective in
this limit. \parn
b) Let us consider, for example, the case $d=2$ and the values
of $n$ in Table 1. For the integer values $n=2,4,7,16,31$
we have determined the analytic expression of $\Kap^{F}_{n 2}$ using Eq. \rref{gv}, and then
maximized this function numerically; the maxima occur, respectively, at $(p, \sigma) \simeq$ $(0.511, 1.05)$,
$(0.417, 0.309)$, $(0.371,0.148)$, $(0.331, 0.0582)$, $(0.316, 0.0290)$. For the large values
$n = 61, 121$, we have used directly the lower bound $K^{FF}_{n 2} = \Kap^{F}_{n 2}(1/(2 \sqrt{2}),
3/(4 n))$; since $n$ is integer, in principle this could be obtained again from Eq.
\rref{gv}, but in these two cases it is more convenient to compute it numerically,
starting from the integral representation \rref{giveby} of $\| f_{p \sigma 2} \|_n$
(note that this contains the non elementary function $I_0$). \parn
For $n = 5/4, 3/2, 5/2$, $\Kap^{F}_{n 2}(p, \sigma)$
has been computed numerically for many sample values of $(p, \sigma)$,
starting again from \rref{giveby}; in this case, approximate
maximization has been performed choosing the best value in the sample. The maxima are
attained at $p \simeq 0.354$ in the three cases, and $\sigma \simeq 5.22, 2.41, 0.696$, respectively. \parn
The numerical computation of $\Kap^{F}_{n 2}$ and $K^{F}_{n d}$ is difficult for the small values $n = 1 + 10^{-4}$,
$1 + 10^{-3}$ and $1 + 10^{-1}$. On the other hand, for the reasons already explained at the end
of the previous section the Fourier bounds should be below the Bessel bounds for these extreme
values of $n$; therefore to construct Table 1 in these cases we have given up computing $K^{F}_{n 2}$,
and we have chosen directly for $K^{-}_{n 2}$ a Bessel lower bound.
\vskip 0.2cm \noindent
\appendix
\section{Appendix. The integral $\boma{I_{\nnu \mmu}(h)}$.}
\label{appint}
This integral is defined by Eq. \rref{imunu}; we want to prove Eq. \rref{result}. We start
from the identity (\cite{Wat}, page 440)
\beq K^2_{\mmu/2}(r) = 2 \int_{0}^{+\infty} d t \, K_{\mmu}(2 r \cosh t)~, \feq
and insert it into \rref{imunu}; this gives
\beq I_{\nnu \mmu}(h) = 2 \int_{0}^{+\infty} d t \, \int_{0}^{+\infty} d r \, r^{\nnu + \mmu + 1} J_{\nnu}(h r)
K_{\mmu}(2 r \cosh t)~. \label{bfw} \feq
On the other hand (\cite{Wat}, page 410)
\beq \int_{0}^{+\infty} d r \, r^{\nnu + \mmu + 1} J_{\nnu}(h r)
K_{\mmu}(2 r \cosh t) =
{\Gamma(\nnu + \mmu + 1) h^{\nnu} \over 2^{\nnu + 2} \cosh^{2 \nnu + \mmu + 2} t} \times \feq
$$ \times F\Big(\nnu + \mmu + 1,
\nnu + 1, \nnu + 1; - {h^2 \over 4 \cosh^2 t}\Big)
= {\Gamma(\nnu + \mmu + 1) h^{\nnu} \over 2^{\nnu + 2} \cosh^{2 \nnu + \mmu + 2} t}
\left( 1 + {h^2 \over 4 \cosh^2 t} \right)^{-\nnu - \mmu - 1}~, $$
where the last passage depends on \rref{ey}. Returning to Eq. \rref{bfw} we obtain
$$ I_{\nnu \mmu}(h) = {\Gamma(\nnu + \mmu + 1) h^{\nnu} \over 2^{\nnu + 1}} \int_{0}^{+\infty} d t \,
{1 \over \cosh^{2 \nnu + \mmu + 2} t} \left( 1 + {h^2 \over 4 \cosh^2 t} \right)^{-\nnu - \mmu - 1} = $$
\beq = {\Gamma(\nnu + \mmu + 1) h^{\nnu} \over 2^{\nnu + 2}} \int_{0}^{1} d s \,
s^{\nnu + \mmu/2} (1 - s)^{-1/2} ( 1 + {h^2 \over 4} s )^{-\nnu - \mmu - 1}~, \feq
the last passage following with the change of variable $s = 1/\cosh^{2} t$. Now, comparison with
\rref{irep} gives the thesis \rref{result}.
\vskip 0.2cm \noindent
\section{Appendix. Proof of Prop. \ref{mainp} on Laplace integrals.}
\label{appeaa}
We recall the notations and assumptions (\ref{int}-\ref{supp}), and point out some consequences
of our hypotheses.
\parn First of all, by the monotonicity of $\ffi$,
$\ffi(b):= \lim_{t \vain b^{-}} \ffi(t)$ exists in $(0, +\infty]$,
and $\ffi$ is a $C^1$ diffeomorphism between $(0,b)$ and $(0, \ffi(b))$. \parn
Moreover, by Eq. \rref{supp}, there are a constant $\ep \in (0,b)$ and a bounded function
$\beta \in C((0,\ep), \reali)$ such that
\beq \xi(t) = \sum_{i=0}^{\ell-1} P_i \, \ffi(t)^{\alpha_i-1} + \beta(t) \ffi(t)^{\alpha_{\ell} -1}
\qquad \mbox{for all $t \in (0,\ep)$}~. \label{bet} \feq
Putting the attention to Eq. \rref{int} and dividing integration in two parts,
we get
\beq L(n) = \II(n) + \JJ(n)~, \label{inte} \feq
$$ \II(n) :=
\int_{0}^{\ep} d t~\ppsi(t)~ e^{-n \ffi(t)}~, \qquad \JJ(n) :=
\int_{\ep}^{b} d t~\ppsi(t)~ e^{-n \ffi(t)}~. $$
Let us estimate $\II(n)$. Introducing the new variable $s = \ffi(t)$ and then using \rref{bet} we obtain
\beq \II(n) = \int_{0}^{\ffi(\ep)} d s
~\xi(\ffi^{-1}(s))~ e^{-n s} = \sum_{i=0}^{\ell-1} P_i \II_i(n) + \delta \II_\ell(n)~, \feq
$$ \II_i(n) := \int_{0}^{\ffi(\ep)} d s~s^{\alpha_i-1} e^{-n s}~, \qquad \delta \II_{\ell}(n) :=
\int_{0}^{\ffi(\ep)} d s~s^{\alpha_\ell-1} \beta(\ffi^{-1}(s)) e^{-n s}~. $$
The above integrals are related to the incomplete Gamma function
\beq \gamma(\alpha, u) := \int_{0}^u d v\, v^{\alpha-1} e^{-v} =
\Gamma(\alpha) - \int_{u}^{+\infty} d v\, v^{\alpha-1} e^{-v} = \feq
$$ = \Gamma(\alpha) + O\left( u^{\alpha-1} e^{-u} \right) \qquad
\mbox{for $u \vain +\infty$} \qquad (\alpha > 0)~
$$
(concerning the asymptotics of $\gamma$ for $u \vain +\infty$, see \cite{Olv}). As for $M_i$, with a
variable change $s = v/n$ we get
\beq \II_i(n) = {\gamma(\alpha_i, n \ffi(\ep)) \over n^{\alpha_i}} =
{\Gamma(\alpha_i) \over n^{\alpha_i}} + O\left({e^{-n \ffi(\ep)} \over n} \right)
\quad \mbox{for $n \vain +\infty$} ~; \label{oi} \feq
furthermore,
\beq | \delta \II_{\ell}(n) | \leqs (\sup_{(0, \ep)} | \beta |) ~
\int_{0}^{\ffi(\ep)} d s~s^{\alpha_\ell-1} e^{-n s}
= (\sup_{(0, \ep)} | \beta |) ~{\gamma(\alpha_{\ell}, n \ffi(\ep)) \over n^{\alpha_\ell}} =
\label{ol} \feq
$$ = (\sup_{(0, \ep)} | \beta |) ~
\left[{\Gamma(\alpha_\ell) \over n^{\alpha_\ell}} + O\left({e^{-n \ffi(\ep)} \over n} \right)\right]=
O\left({1 \over n^{\alpha_\ell}}\right) \quad \mbox{for $n \vain +\infty$}~. $$
To estimate $\JJ(n)$, we fix $n_1 > n_0$ and write
$\JJ(n) = \int_{\ep}^{b} d t~\ppsi(t)~ e^{-(n - n_1) \ffi(t)} e^{-n_1 \ffi(t)}$;
for all $n \in [n_1, +\infty)$, this implies
\beq | \JJ(n) | \leqs
e^{-(n - n_1) \ffi(\ep)} \int_{\ep}^{b} d t~| \ppsi(t) |~ e^{-n_1 \ffi(t)}
= O\left(e^{-n \ffi(\ep)}\right) \quad \mbox{for $n \vain +\infty$}  \label{oj} \feq
(recall that $0 < \ffi(\ep) \leqs \ffi(t)$ for $t \in [\ep, b)$).\parn
From Eq.s \rref{inte}, \rref{oi}, \rref{ol} and \rref{oj} we get the thesis \rref{maineq}. \fine
\vskip 0.4cm \noindent
\textbf{Acknowledgments.} We acknowledge D. Bambusi and S. Paveri Fontana
for useful bibliographical indications.
This work has been partially supported by the GNFM
of Istituto Nazionale di Alta Ma\-te\-ma\-ti\-ca and by MIUR,
Research Project Cofin/2004 "Metodi geometrici nella teoria delle onde non lineari e applicazioni".

\end{document}